\documentclass[opre,nonblindrev]{informs3_nojournal} 

\OneAndAHalfSpacedXII 

\usepackage{endnotes}
\let\footnote=\endnote

\usepackage{ulem}
\usepackage{natbib}
 \bibpunct[, ]{(}{)}{,}{a}{}{,}%
 \def\bibfont{\small}%
\TheoremsNumberedThrough     
\ECRepeatTheorems

\EquationsNumberedThrough    

\usepackage[english]{babel}
\usepackage[utf8]{inputenc}
\usepackage{algorithm}
\usepackage[noend]{algpseudocode}
\usepackage{mathrsfs}
\usepackage{soul}
\usepackage{bm}
\usepackage{hyperref}	

\usepackage{makecell} 
\usepackage{array}
\newcolumntype{C}[1]{>{\centering\arraybackslash}m{#1}} 

\begin{document}


\RUNAUTHOR{}

\RUNTITLE{Data-Driven Bayesian Nonparametric Approach for Black-Box Optimization}

\TITLE{A Data-Driven Bayesian Nonparametric Approach for Black-Box Optimization}

\ARTICLEAUTHORS{%
	\AUTHOR{Haowei Wang}
	
	\AFF{Rice-Rick Digitalization PTE. Ltd., 51 Goldhill Plaza, Singapore 308900, SINGAPORE \EMAIL{wanghaowei58@gmail.com}}
	
	\AUTHOR{Xun Zhang}
	
	\AFF{School of Business, National University of Singapore, Singapore 117576, SINGAPORE \EMAIL{e0225088@u.nus.edu}}
	
	\AUTHOR{Szu Hui Ng}
	
	\AFF{ Department of Industrial Systems Engineering and Management,  National University of Singapore, Singapore 117576, SINGAPORE \EMAIL{isensh@nus.edu.sg}}
	
	\AUTHOR{Songhao Wang}
	
	\AFF{SUSTech Business School, Southern University of Science and Technology, 1088 Xueyuan Avenue, Shenzhen 518055, CHINA \EMAIL{wangsh2021@sustech.edu.cn}}
	
}

\ABSTRACT{%
We present a data-driven Bayesian nonparametric approach for global optimization (DaBNO) of stochastic black-box function. The function value depends on the distribution of a random vector. However, this distribution is usually complex and hardly known in practice, and is often inferred from data (realizations of random vectors). 
The DaBNO accounts for the finite-data error that arises when estimating the distribution and relaxes the commonly-used parametric assumption to reduce the distribution-misspecified error.
We show that the DaBNO objective formulation can converge to the true objective asymptotically. 
We further develop a surrogate-assisted algorithm DaBNO-K to efficiently optimize the proposed objective function based on a carefully designed kernel. Numerical experiments are conducted with several synthetic and practical problems, demonstrating the empirical global convergence of this algorithm and its finite-sample performance.}%


\maketitle

%

\section{Introduction} \label{sec3:introduction}
We focus on stochastic global optimization for expensive stochastic black-box functions:
\begin{equation}
	\min_{\BFx \in \mathcal{X}}  \quad  f(\BFx, \BFP^c) \triangleq \text{E}_{\BFu \sim \BFP^c}[h(\BFx,\BFu)],  \label{trueobjective}
\end{equation}
where $\BFx$ is the decision vector in a compact set $\mathcal{X} \subset \mathbb{R}^d$, $\BFu$ is a random vector (uncontrollable by decision-makers) governed by the probability distribution $\BFP^c$, and $h(\BFx,\BFu)$ is the black-box function value. We are interested in optimizing the expectation of $h(\BFx,\BFu)$ with respect to the distribution of $\BFu$, denoted by $f(\BFx, \BFP^c)$.  Black-box optimization often arises in the design of systems when their performances can be obtained through high-fidelity numerical evaluations, such as the finite element analysis (FEA) for partial differential equations (PDE) in mechanical design, and the simulation models for operational management. Problem \eqref{trueobjective} depicts a class of black-box problems involving random parameters. Two examples stemming from engineering design and stochastic simulation are presented below.

\emph{EXAMPLE 1 (Design of an Acoustic Horn)} The decision makers need to select several design parameters $\BFx$, including geometric parameters and boundary parameters, for an acoustic horn to maximize its expected efficiency $f$ \citep{ng2014multifidelity}. The efficiency is decided as the reflection coefficient $h$, measuring the ratio of how much a wave is internally reflected over that is transmitted out. Here, the random vector $\BFu$ is the uncertain wave number due to the unknown operating condition or the environmental variable when the horn is used in reality. Given $\BFx$ and $\BFu$, $h$ can be evaluated through FEA solving the governing PDEs.

\emph{EXAMPLE 2 (Simulation of an Inventory System)} The decision makers use stochastic simulations to study an inventory system. They need to decide the ordering level $\BFx$ to maximize the expected profit $f$. Here, the random vector $\BFu$ is the customer demand. Given $\BFx$ and $\BFu$, they can run the simulation to mimic the operations of the systems to obtain profit $h$. For stochastic simulations, the random vector is also called `input parameter' to the model \citep{zhou2017simulation}.

\subsection{Motivation}
One challenge in stochastic optimization is that solving \eqref{trueobjective} directly is often impossible in practice, as $\BFP^c$ is hardly known exactly. We can instead solve some approximations to \eqref{trueobjective} based on finite realizations of $\BFu$. 
The recent Bayesian framework \citep{wu2018bayesian} provides an asymptotically valid approximation to adequately account for the finite-data uncertainty in estimating $\BFP^c$, where the authors  assume a parametric form for $\BFP^c$ with parameter $\BFlambda^c$. The posterior of $\BFlambda^c$, denoted as $\BFpi^c$, quantifies the estimation uncertainty. They further proposed to solve the following problem instead of \eqref{trueobjective}:
\begin{equation}
	\min_{\BFx\in \mathcal{X}} \quad  \text{E}_{\BFlambda \sim \BFpi^c} \left[ \text{E}_{\BFu\sim \BFP_{\BFlambda}} [h(\BFx,\BFu)]\right].                               \label{objectivefunction}
\end{equation}

The existing research applying such a Bayesian approach assume a pre-defined parametric form for the distribution of $\BFu$. This however can be too restrictive in practice with limited prior information. Also, real world data usually exhibit rich properties such as multi-modality, skewness, and excess kurtosis \citep{xie2019bayesian}. These diverse and complex natures make it difficult to accurately represent the distribution of $\BFu$ with standard parametric families. It can lead to model misspecification errors when parametric distributions are assumed, even when the size of real world dataset is large.  

Another challenge lies in the high evaluation cost for the black-box functions. Due to the complex nature of the underlying real-world systems, it is not uncommon to  take several hours or days to  run some of these experiments. Additionally, the black-box function can be multimodal with no gradient information. It is thus of interest to develop efficient algorithms that can find the optimal solutions with limited budget which have global convergence guarantees.  

Motivated by limited real world data and the complex natures of the uncertain distributions, we develop a non-parametric global optimization approach for stochastic black-box systems. Specifically, we design an approximation \eqref{objfunction_0} to the true objective function \eqref{trueobjective} as follows:
\begin{equation}
	\min_{\BFx\in \mathcal{X}} \text{E}_{\BFP \sim \BFpi} \left[ \text{E}_{\BFu \sim \BFP} [h(\BFx,\BFu)] \right].   \label{objfunction_0}
\end{equation}
Here, $\BFpi$ is the nonparametric posterior of the distribution $\BFP$. This generalizes \eqref{objectivefunction} by relaxing the parametric assumption. Considering the high cost in black-box evaluations, we further develop an efficient surrogate-assisted algorithm to solve \eqref{objfunction_0}.

\subsection{Literature Review}
Several approximations to problem \eqref{trueobjective} have been proposed based on real-world data. With the empirical distribution, denoted as $\hat{\BFP}$, the Sample Average Approximation (SAA) to \eqref{trueobjective} is defined as \citep{kleywegt2002sample, kim2015guide}:
\begin{equation}
	\min_{\BFx\in \mathcal{X}} \quad \text{E}_{\BFu \sim \hat{\BFP}}[h(\BFx,\BFu)]. \label{approproblem}
\end{equation}
When faced with finite data, the estimation $\hat{\BFP}$ unavoidably contains estimation error. Thus, it is likely that the objective function values of the best decision that solves the problem (\ref{approproblem}) can be quite different from that of (\ref{trueobjective}), resulting in inferior performance. 

One framework that accounts for finite-data uncertainty in estimating $\BFP^c$ is distributionally robust optimization (DRO). DRO selects an ambiguity set of distributions, denoted as $D$, which can be built by choosing distributions from a ball centered around a reference distribution.
The following approximated objective function is then built to select $\BFx$ considering the worst case in $D$:
\begin{equation}
	\min_{\BFx\in \mathcal{X}} \max_{{\BFP} \in D} \quad \text{E}_{\BFu \sim {\BFP}}[h(\BFx,\BFu)]. \label{dro}
\end{equation}
DRO has been extensively studied for ordinary mathematical programming where the objective functions have explicit formulations \citep{bertsimas2018data,goh2010distributionally,mohajerin2018data,wiesemann2014distributionally} (in contrast to the black-box setting in this work) since the work of \cite{scarf1957min}. We refer to \cite{rahimian2019distributionally} for a thorough review. Recently, \cite{kirschner2020distributionally} and \cite{husain2022distributionally} propose surrogate-based algorithms to solve black-box optimization with ambiguity sets based on the maximum mean discrepancy and $\phi$-divergence, respectively. One possible issue of DRO is that the solutions can be conservative if the ambiguity set is not properly defined \citep{wang2016likelihood}, as the performance under true uncertain distributions may be much better than that of the worst case in $D$.

The Bayesian approach \eqref{objectivefunction} is proposed by \cite{zhou2015simulation}. They empirically showed that the solution to \eqref{objectivefunction} outperforms that of \eqref{approproblem} when the size of the real world observed dataset { is small and $h$ is sensitive to $\lambda$. \cite{wu2018bayesian} further prove the consistency and asymptotic normality of the optimal value} of \eqref{objectivefunction} when datasize tends to infinity. They further develop Bayesian risk optimization to explores a balance between an over-optimistic solution of \eqref{approproblem} (totally ignoring the uncertainty in $\hat{\BFP}$) and the more conservative DRO solution of \eqref{dro}. { \cite{Pearce2018} and \cite{wang2019} developed Kriging surrogate-based algorithms to solve \eqref{objectivefunction}. }

The current Bayesian approaches mostly assume that the distribution family of the random parameter $\BFu$ is known, which may not be reliable in practice due to the complexity of the real world process. To mitigate the distribution selection error,  the Bayesian Model Averaging (BMA) \citep{chick2001input} has been proposed to quantify the distribution family uncertainty by maintaining a posterior average over a set of distribution family choices. However, it can still be difficult to select the candidate parametric families in the first place. Others alleviate this through nonparametric modeling including resampling, bootstrapping \citep{barton2001resampling,barton2013quantifying} and Dirichlet process mixtures \citep{xie2019bayesian}. The main focus of these works, however, is to quantify the finite-data uncertainty of a single experiment instead of  \textbf{\emph{global optimization}}. Thus, the question remains on how to solve black-box optimization problems with unknown and complex distributions of $\BFu$. 

When it comes to numerically optimize expensive black-box functions such as \eqref{objectivefunction} and  \eqref{approproblem},  surrogate-assisted optimization algorithms are commonly used  \citep{forrester2008engineering,jones2001taxonomy}, including algorithms based on response surface methodology \citep{myers2016response}, radial basis functions \citep{regis2007stochastic,muller2017socemo,muller2019surrogate}, and the Kriging model  \citep{jones1998efficient,meng2022combined}. The surrogate model is a statistical model approximating the true expensive function which can be used to guide the search for the optimum. 
The surrogate-based algorithm typically operates in an iterative manner. In each iteration, a new point $\BFx$ is selected for expensive evaluation based on the current surrogate model, which will be updated when the evaluation results are returned. The evaluation point $\BFx$ can be chosen through an inner optimization problem, such as optimizing the surrogate or an improvement function \citep{jones1998efficient}. For black box function $h(\BFx, \BFu)$, we need to specify not only $\BFx$ but $\BFu$ to run the black-box function. This requires, in each iteration, additionally selecting $\BFu$ or a distribution of $\BFu$ (so that realizations of $\BFu$ can be generated). In this paper, we propose a Kriging-based optimization algorithm which takes into account the uncertainty while has theoretical guarantee to solve \eqref{objfunction_0}. Readers can refer to  \cite{shahriari2016taking} for a thorough review on Kriging-based optimization algorithm.  

\subsection{Our Contributions}
In this work, we propose a Data-driven Bayesian Nonparametric Optimization (DaBNO) approach for stochastic black-box systems. The key contributions of our paper can be summarized as follows.

1. We specifically design a \textbf{\emph{nonparametric approximation}} \eqref{objfunction_0} to the true objective function \eqref{trueobjective} based on the Dirichlet process. 
This extends the parametric approach taken by \cite{wu2018bayesian} to enable more flexible capturing of the distribution uncertainty.
We show the \textbf{\emph{consistency}} of this approximation: the objective function, as well as the optimal solution, converges to the true objective, as the number of data points ${\BFu}$ tends to infinity. These convergence results generally ensure that solving the approximations is equivalent to solving \eqref{trueobjective} as the size of the data set increases. 
Moreover, in an online setting where additional data of ${\BFu}$ are provided sequentially during the decision-making process \citep{wu2022data}, a consistent objective function corrects and guides the optimization towards the right direction. 

2. We show the \textbf{\emph{asymptotic normality}} of the objective function, optimal solution, and optimal value. Asymptotic normality provides additional insights on how quick the convergence is as well as the distribution of the solutions. Specifically, we show that $\sqrt{s}$ ($s$ is the number of observations of $\BFu$) times the difference between our estimated optimal value/solution and the true ones under perfect uncertain distribution follows a normal distribution asymptotically as the real world data size increases. This reveals a rate of convergence at $O_P(1/\sqrt{s})$. Asymptotic normality results can be used for statistical inferences (e.g., hypothesis testing and building confidence intervals) to the optimal solution and objective \eqref{trueobjective}. 
This nonparametric objective function, as well as its asymptotic properties, is potentially of independent interest in mathematical programming for non-black-box functions. 

3. We develop a \textbf{\emph{surrogate-assisted algorithm}} based on the Kriging model to practically solve \eqref{objfunction_0}.  Specifically, we build a Kriging model of the function $f(\BFx,\BFP^c)$ with respect to $\BFx$ and the random vector distribution $\BFP$ based on a mix of Euclidean distance and Wassertain distance. Defining \eqref{objfunction_0}  as the DaBNO objective function, a model of it is obtained through aggregating $\BFP$ according to its posterior. Finally, we design a searching criterion to jointly select the next $\BFx$ and $\BFP$, that  have the largest improvement in the DaBNO objective. We also prove the convergence of this algorithm. This algorithm hence provides a practical approach to solving problems in  engineering design, operations research, and computer science with expensive objectives and uncertain inputs.


A preliminary version of this paper was published in \cite{wang2020gaussian}. We significantly improve the conference paper both theoretically and empirically. This manuscript formalizes and extends the framework to general black-box systems which was only superficially presented in \cite{wang2020gaussian}. In addition, we provide the convergence proofs and asymptotic properties of the proposed approximated objective function in Section \ref{sec3:NBRO}. The consistency property of the surrogate-assisted algorithm based on the Kriging model has been substantiated through a detailed proof in Section \ref{convergence}. A more comprehensive numerical study is provided to illustrate the work.

The remaining structure of this paper is as follows. Section \ref{sec3:NBRO} formulates the objective function in DaBNO and its asymptotics. Section \ref{sec3:algo} develops the Kriging surrogate model for DaBNO. Section \ref{EGO} presents the new algorithm DaBNO-K for solving DaBNO. Section \ref{convergence} proves the theoretical convergence properties of the DaBNO-K algorithms. In Section \ref{sec3:numericalexperiments}, we empirically test the framework and algorithm with several synthetic and practical problems. Section \ref{sec3:conclusion} concludes the paper.

\section{Objective function for DaBNO and its asymptotics} \label{sec3:NBRO}
We formulate the DaBNO objective and its non-parametric modeling with Dirichlet Process in Section \ref{sec3_1}. The Dirichlet Process is a flexible approach for Bayesian nonparametric modeling with efficient updating and sampling mechanisms. Section \ref{sec3:properties} presents the asymptotics of this objective function.

\subsection{The DaBNO objective function}
\label{sec3_1}
{To hedge the risk of solving Equation (\ref{approproblem}) based on $\hat{\BFP}$, we model $\BFP^c$ with a posterior $\BFP \sim \BFpi$ and propose to solve
	\begin{equation}
		\min_{\BFx\in \mathcal{X}} \quad g(\BFx)=\text{E}_{\BFP \sim \BFpi} \left[ \text{E}_{\BFu \sim \BFP} [h(\BFx,\BFu)] \right].   \label{objfunction}
	\end{equation}
	The inner expectation in (\ref{objfunction}) is taken with respect to distribution of $\BFu$ and the outer expectation is taken against the posterior distribution of $\BFP$} to capture the finite data uncertainty. We note here that if $\BFP$ is chosen from a parametric family, the DaBNO formulation reduces to \eqref{objectivefunction}.

Denote the real world data set by $\BFD_s= \{\BFu_1, \cdots, \BFu_s\}$ where the data points are identically and independently distributed with $\BFP^c$. From a Bayesian perspective, the unknown $\BFP^c$ is modeled as a random distribution, denoted by $\BFP$, with a prior $\BFpi_0(\BFP)$. This prior represents our initial belief and knowledge about $\BFP^c$. 
It will be updated after observing the data set $\BFD_s$ to obtain the posterior $\BFpi(\BFP|\BFD_s)$. If no initial knowledge is available, non-informative priors can be used.
In this paper, we assume non-informative priors and propose to use the conjugate Dirichlet process prior \citep{ferguson1973bayesian}, a popular prior for nonparametric modeling. 
The Dirichlet process (DP) is a family of stochastic processes whose realizations are probability distributions. It can be denoted as $DP(\alpha,\BFP^0)$ where {$\alpha$ is the concentration parameter and $\BFP^0$ is the base distribution}. 
The base distribution $\BFP^0$ {can be thought of as a guess for $\BFP^c$.} It is the mean value of the Dirichlet prior, which draws distributions ``around'' it.  
{High values of $\alpha$ imply high confidence in $\BFP^0$}  and the prior is more concentrated around $\BFP^0$.  
Given  $\BFu_1,...,\BFu_s\sim \BFP^c$ and the prior $DP(\alpha,\BFP^0)$,  the posterior is also a Dirichlet process: \[\BFpi(\BFP|\BFu_1,...,\BFu_s) \sim DP(\alpha+s, \frac{\alpha \BFP^0 + \sum_{j=1}^{s}\delta_{\BFu_j}}{\alpha+s}). \]
where $\delta$ is the Dirac delta measure
\begin{equation*}
	\delta_{\BFu_i}(E)=\left\{\begin{array}{ll}
		1 & \text { if } \BFu_i \in E \\
		0 & \text { if } \BFu_i \notin E
	\end{array}\right. ,
\end{equation*}
where $E$ is any event in the event space. The base distribution of the posterior Dirichlet process  becomes $\frac{\alpha \BFP^0 + \sum_{j=1}^{s}\delta_{\BFu_j}}{\alpha+s} = \frac{\alpha}{\alpha+s} \BFP^0+ \frac{s}{\alpha+s}(\frac{1}{s}\sum_{j=1}^{s}\delta_{\BFu_j})$ which is a weighted average of the prior base distribution $\BFP^0$ and the empirical distribution $\frac{1}{s}\sum_{j=1}^{s}\delta_{\BFu_j}$. DP has been successfully applied to capture rich natures of distributions in a wide range of applications, such as modeling the distribution of patient failure time \citep{christensen1988modelling}, molecular lineage-specific rates of substitution \citep{heath2012dirichlet} and the triglyceride levels of patients \citep{barcella2018dependent}.  It is easy and cheap to obtain samples from the Dirichlet Process posterior using the algorithm of stick breaking process \citep{hjort2010bayesian}.
The formal definition and construction of a Dirichlet process is based on the language of random measures which we omit for ease of exposition. 
The readers may check Chapter 4 of \cite{ghosal2017fundamentals}) for detailed descriptions.


\subsection{Asymptotic properties of the DaBNO objective function}\label{sec3:properties}
In this section, we study the consistency and asymptotic normality results of the objective function in DaBNO. For consistency, we show the following three \emph{almost surely} results as the size $s$ (the number of real world observations $\BFu$) tends to infinity: 
\begin{itemize}
	\item [1)] the objective function of DaBNO, $\mathbb{E}_{\BFP \sim \BFpi}[\mathbb{E}_{\BFu \sim \BFP} [h(\BFx,\BFu)]]$, converges pointwisely to  true objective function $\mathbb{E}_{\BFu \sim \BFP^c}[h(\BFx,\BFu)]$;
	\item [2)] the optimal solution of DaBNO objective function, $\BFx_s=\arg\min_{\BFx}\mathbb{E}_{\BFP \sim \BFpi}[\mathbb{E}_{\BFu \sim \BFP} [h(\BFx,\BFu)]]$, converges to the true optimal design: $\BFx^*=\arg\min_{\BFx}\mathbb{E}_{\BFu \sim \BFP^c}[h(\BFx,\BFu)]$; 
	\item [3)] the optimal value of DaBNO objective function, $\mathbb{E}_{\BFP \sim \BFpi}[\mathbb{E}_{\BFu \sim \BFP} [h(\BFx_s,\BFu)]]$, converges to the true optimal value $\mathbb{E}_{\BFu \sim \BFP^c}[h(\BFx^*,\BFu)]$.
\end{itemize}
These consistency results show that the DaBNO problem in (\ref{objfunction}) will approach the true one in (\ref{trueobjective}) as $s$ gets large, even if the true underlying distribution $\BFP^c$ is unknown. In the ideal case of a sufficiently large sample size of the real observations, solving the DaBNO problem is asymptotically equivalent to solving the true problem. 

In addition to the consistency, we further investigate the asymptotic normality property of the solutions as well as the function value of \eqref{objfunction}. Asymptotic normality results describe the distribution of the solutions to the problem as the number of real world observations increases. They can be used to examine how quick the convergence is and facilitate statistical inferences. For instance, we can build confidence intervals for the optimal value/solution to quantify the uncertainty in using DaBNO. Furthermore, one can specify a tolerance level of the variance (or the width of the confidence interval), based on which they can determine the minimum amount of real world data needed to achieve the target variance with the estimators of the asymptotic variances. In this section, we show that $\sqrt{s}$ times the difference between our estimated optimal value/solution and the true ones under perfect uncertain distribution follows a normal distribution asymptotically, indicating a rate of convergence at $O_P(1/\sqrt{s})$. To the best of our knowledge, the only piece of work studying the asymptotic normality of the estimators in the literature of input uncertainty is \cite{wu2018bayesian}. However, as discussed later, \cite{wu2018bayesian} only studies parametric model, and their techniques do not extend here.

Both consistency and asymptotic normality are theoretical evidence supporting the assertion that DaBNO objective is a well founded approximation to the original objective function when the uncertain distribution is modeled using the Dirichlet process. Before proceeding to the specific theorems and proofs, we make the following assumption about $h(\BFx,\BFu)$ and $\BFu$.

%

\begin{assumption} \label{assum}
	\item [ \quad (i)] $h(\BFx,\BFu)$ is Lipschitz in $\BFx$ and $\BFu$, i.e., there exists positive constants $L_1$, $L_2$, and $L_3$ such that 
	\begin{equation*}
		\begin{aligned}
			|h(\BFx,\BFu_1)-h(\BFx,\BFu_2)| &\leq L_1||\BFu_1-\BFu_2||, \quad \forall \BFx\\
			| h(\BFx_1,\BFu)- h(\BFx_2,\BFu)|&\leq L_2||\BFx_1-\BFx_2||,\quad \forall \BFu.\\
			||\partial_{\BFx} h(\BFx_1,\BFu)-\partial_{\BFx} h(\BFx_2,\BFu)||&\leq L_3||\BFx_1-\BFx_2||,\quad \forall \BFu.
		\end{aligned}
	\end{equation*}
	where $\partial_{\BFx}$ means taking the partial derivative with respect to $\BFx$ and $\| \cdot \|$ denotes the  Euclidean norm. 
	\label{assum:1}
	\item [ \quad (ii)] $\BFu$ has bounded support, i.e., there exists a compact set $\Omega=[\underline{u}_1,\overline{u}_1]\times [\underline{u}_2,\overline{u}_2] \times\cdots \times [\underline{u}_l,\overline{u}_l]$ such that $\BFu\in \Omega$, where $l$ is the dimension of $\BFu.$ 
	\label{assum:2}
\end{assumption}
The condition on the Lipschitz continuity limits how fast the function $h(\BFx,\BFu)$ can change with respect to both $\BFx$ and $\BFu$, which is an assumption often imposed in black-box models \citep{fan2018surrogate,ghosh2019robust}. The assumption of bounded support for $\BFu$ also holds in many situations. Take an inventory problem as an example, the random demand is usually bounded.

To establish the consistency and asymptotic normality results, we first transform the objective function in \eqref{objfunction} into a summation of several independent and identically distributed (i.i.d.) terms in Lemma \ref{lemma:1}. Based on this representation, we analyze the asymptotic behavior from the viewpoint of empirical process. In this paper, all the theoretical proofs are provided in the Appendix.
\begin{lemma}\label{lemma:1}
	Suppose Assumption \ref{assum} holds. Then
	\begin{equation*}
		\mathbb{E}_{\BFpi}[\mathbb{E}_{\BFP} [h(\BFx,\BFu)]]=\mathbb{E}_{\BFP^0}\frac{\alpha}{\alpha+s}h(\BFx,\BFu)+\frac{1}{\alpha+s}\sum_{j=1}^s  h(\BFx,\BFu_j).  
	\end{equation*}   
\end{lemma}
Now we can prove the consistency and asymptotic normality given by the following two theorems. 
Suppose that $\BFx_s=\arg\min_{\BFx}\mathbb{E}_{\BFP \sim \BFpi}[\mathbb{E}_{\BFP}[ h(\BFx,\BFu)]] $ and $\BFx^*=\arg\min_{\BFx} \mathbb{E}_{\BFP^c} [h(\BFx,\BFu)]$. The following  Theorem \ref{theorem:1} and   Theorem \ref{theorem:2} summarize the results of consistency and the asymptotic normality respectively. We use $\to_{a.s.}$ to denote the almost sure convergence and $\leadsto$ to denote the convergence in distribution. 

\begin{theorem}\label{theorem:1}
	Suppose Assumption \ref{assum} holds. We have
	\begin{itemize}
		\item[(i)] (Consistency of the value of the objective function) As $s\to+\infty$,
		\begin{equation*}
			\mathbb{E}_{\BFP \sim \BFpi}[\mathbb{E}_{\BFP} [h(\BFx,\BFu)]] \to_{a.s.} \mathbb{E}_{\BFP^c}[h(\BFx,\BFu)] .
		\end{equation*}
		\item[(ii)] (Consistency of the optimal solution) 
		As $s\to +\infty$,
		\begin{equation*}
			\BFx_s \to_{a.s.} \BFx^*.
		\end{equation*}
		\item[(iii)] (Consistency of the optimal value) As $s\to+\infty$,
		\begin{equation*}
			\mathbb{E}_{\BFP \sim \BFpi}[\mathbb{E}_{\BFP} [h(\BFx_s,\BFu)]]\to_{a.s.} \mathbb{E}_{\BFP^c}[h(\BFx^*,\BFu)].
		\end{equation*}
	\end{itemize}
\end{theorem}

Denote the derivative matrix of ${\BFx}\mapsto \mathbb{E}_{\BFu\sim {\BFP}_c}\frac{\partial}{\partial {\BFx}}f({\BFx},\BFu)$ at ${\BFx}^*$ as $V_{{\BFx}^*}$.
We summarize the asymptotic normality results below.

\begin{theorem}\label{theorem:2}
	Suppose Assumption \ref{assum} holds, and further assume  that  $V_{{\BFx}^*}$ is invertible. Then, we have
	\begin{itemize}
		\item[(i)] (Asymptotic normality of the value of the objective function.) As $s\to+\infty$,
		\begin{equation*}
			\sqrt{s}\left(\mathbb{E}_{\BFP\sim \BFpi}[\mathbb{E}_{{{\BFP}}} [h({{\BFx}},{\BFu})]]-\mathbb{E}_{{{\BFP}}^c}[h({{\BFx}},{\BFu})]\right) \leadsto N(0,\sigma_{h({\BFx})}^2),
		\end{equation*}
		where $\sigma_{h({\BFx})}^2=\textnormal{Var}(h({{\BFx}},{\BFu}))$.
		\item[(ii)] 	(Asymptotic normality of the optimal solution)
		$$\sqrt{s}({{\BFx}}_s-{{\BFx}}^*)\leadsto N(0,\Sigma)$$
		where $\Sigma=V_{{\BFx}^*}^{-1}\mathbb{E}\frac{\partial}{\partial {\BFx}}h({\BFx^*},\BFu)\frac{\partial}{\partial {\BFx}}h({\BFx^*},\BFu)^T (V_{{\BFx}^*}^{-1})^{T}.$ 
		\item[(iii)]  (Asymptotic normality of the optimal value)  As $s\to+\infty$,
		\begin{equation*}
			\sqrt{s}\left(\mathbb{E}_{\BFP\sim \BFpi}[\mathbb{E}_{{{\BFP}}} [h({{\BFx}}_s,{\BFu})]]-\mathbb{E}_{{{\BFP}}^c}[h({{\BFx}}^*,{\BFu})]\right) \leadsto N(0,\sigma_{h({\BFx}^*)}^2),
		\end{equation*}
		where $\sigma_{h({\BFx}^*)}^2=\textnormal{Var}(h({{\BFx}}^*,{\BFu})).$ 
	\end{itemize}
	
\end{theorem}

\begin{remark} 
	Our results coincide with the results in \cite{wu2018bayesian} for a special case. Here, we consider an example when $\BFu$ is a binary variable with $P(\BFu=\BFu_1)=\theta^c$ and $P(\BFu=\BFu_2)=1-\theta^c$.
	According to our Theorem \ref{theorem:2}, we have
	$
	\sqrt{s}\left(\mathbb{E}_{\BFpi}[\mathbb{E}_{\BFP}[ h(\BFx,\BFu)]]-\mathbb{E}_{\BFP^c}[h(\BFx,\BFu)]\right) \leadsto N(0,\sigma_{h(x)}^2),
	$
	where $\sigma_{h(x)}^2=\text{Var}(h(\BFx,\BFu)) = \left( \theta^c h^2(\BFx,\BFu_1)+(1-\theta^c)h^2(\BFx,\BFu_2) \right)-(\theta^c h(\BFx,\BFu_1)+(1-\theta^c )h(\BFx,\BFu_2))^2=\theta^c (1-\theta^c )\left( h(\BFx,\BFu_1)-h(\BFx,\BFu_2) \right)^2.$ 
	Since $\BFu$ is a binary random variable with parameter $\theta_c$, we apply  Theorem 4.4 of \cite{wu2018bayesian} to obtain
	$\sqrt{s}\left(\mathbb{E}_{\BFpi}[\mathbb{E}_{\BFP}[ h(\BFx,\BFu)]]-\mathbb{E}_{\BFP^c}[h(\BFx,\BFu)]\right) \leadsto N(0,\sigma^2)$,
	where $\sigma^2=\left(\partial_{\theta} \mathbb{E}h(\BFx,\BFu)\right)^2I(\theta^c)^{-1}=\left( \theta^c  h(\BFx,\BFu_1)+(1-\theta^c )h(\BFx,\BFu_2) \right)'^{2}_{\theta} I(\theta^c )^{-1}=(h(\BFx,\BFu_1)-h(\BFx,\BFu_2))^2I(\theta^c )^{-1} = \theta^c (1-\theta^c )\left( h(\BFx,\BFu_1)-h(\BFx,\BFu_2) \right)^2$ as $ I(\theta^c )=\frac{1}{\theta^c (1-\theta^c )}$ for binary variable $\BFu$. Hence, for this special case, our results for general nonparametric models recover the results of \cite{wu2018bayesian}. Furthermore, when $\BFx=\BFx^*,$ we can also recover the asymptotic variance provided in Theorem 4.13 of \cite{wu2018bayesian}.
	However, we won't have this equality for more general distribution of $\BFu$.
	The asymptotic variance given in Theorem 4.4 of \cite{wu2018bayesian} matches the Cram{\'e}r-Rao bound \citep{vandervarrt}.
	This means that their method is \emph{efficient} under the assumption that a parametric family is known, and so no other estimating methods will have a lower asymptotic variance.
	In other words,  the asymptotic variances given in the above theorem are greater than or equal to those in \cite{wu2018bayesian}.
	It is interesting to observe that although our method is fully non-parametric, for the special random variable, the derived asymptotic variance can be equal to that derived with the knowledge of the parametric family.
\end{remark}

We notice that \cite{wu2018bayesian} proves the asymptotic normality of the objective value and the optimal solution in a parametric approach. Our results are based on a non-parametric framework and provide an additional asymptotic normality guarantee for the optimal value. Proving asymptotic results in our non-parametric setting is technically different from those in \cite{wu2018bayesian}, which heavily rely on the assumption of the parametric family. This means the posterior distribution can be captured by a finite-dimension vector and \cite{wu2018bayesian} assumes that the objective function admits a Taylor expansion around the true distribution (which is a point in the finite-dimensional euclidean space). In contrast, for our nonparametric setting, the underlying distributions fall inside an infinite-dimensional space. It is hard to verify that a Taylor expansion also exists in this space in general, and thus their analysis does not extend here. Instead, we directly express the objective as an i.i.d sum (asymptotically) and resort to empirical process theory to establish the asymptotic normality results. The additional asymptotic normality for the optimal value proposed in this work can be used to quantify the uncertainties and trustworthiness for the optimal value. It can be used to construct the confidence level and determine the size of real world dataset $s$ when the decision makers have desired level of response value and uncertainty tolerance.

\section{Kriging surrogate model for DaBNO} \label{sec3:algo} 
In practice, the real world observations of $\BFu$ are usually limited. We propose an efficient surrogate-assisted algorithm to solve the black-box optimization problem \eqref{objfunction}. In this section, we develop the Kriging surrogate model for DaBNO. 

Denote $F(\BFx,\BFP)=\text{E}_{\BFu \sim \BFP} [h(\BFx,\BFu)]$, then  (\ref{objfunction}) can be written as follows:
{
	\begin{equation}
		\min_\BFx \quad g(\BFx)=\text{E}_{\BFP \sim \BFpi} [f(\BFx,\BFP)],  \label{equal}
	\end{equation}
}
where $\BFx \in \mathbb{R}^d, \BFu \in \mathbb{R}^l$, and $\BFu \sim \BFP$. Here, $\BFpi$ is a  Dirichlet process posterior for $\BFP = [P^{(1)},\cdots,  P^{(l)}]$, where  $P^{(i)}, i = 1, \cdots, l$ are the distributions for different dimensions of $\BFu$.


We proposed to use the Kriging surrogate model that has been shown to be efficient in modeling the objectives \citep{azzimonti2016quantifying,chen2017flexible,owen2017comparison} and solving optimization problems, especially for expensive problems with limited computing budget \citep{frazier2018bayesian}. { To solve (\ref{equal}), a Kriging model for $g(\BFx)$ is required.} However,  we do not have direct output on $g(\BFx)$ or $f(\BFx,\BFP)$, but only the black-box output $h(\BFx,\BFu)$. Here, $h(\BFx,\BFu)$ can be treated as a random observation of $f(\BFx,\BFP)$ due to the random realization of $\BFu$ following $\BFP$. We propose to first construct a Kriging model for $f(\BFx,\BFP)$ (in Section \ref{SK model}), based on which an approximate model for $g(\BFx)$ is developed in Section \ref{sec:gkriging}. Then we could use this approximate model to develop an EGO algorithm to efficiently solve (\ref{equal}) in Section \ref{EGO}.

\subsection{Stochastic Kriging model of $f(\mathbf{x},\mathbf{P})$ }
\label{SK model}
For ease of exposition in the following section, we denote $y(\BFx,\BFP)$ as a random observation of $f(\BFx,\BFP)$, taking value of $h(\BFx,\BFu)$, where $\BFu$ is a random sample from $\BFP$. We then model $y(\BFx,\BFP)$ as a realization of a random process: 
\begin{equation*}
	y(\BFx,\BFP)= F(\BFx,\BFP)+\epsilon = l(\BFx,\BFP)^{\top} \BFbeta + M(\BFx,\BFP) + \epsilon.
\end{equation*}
%
Here, $F(\BFx,\BFP)$ is a Kriging model (also known as a Gaussian process model) used to model the expectation function $f(\BFx,\BFP)$ and $\epsilon$ is a normal noise to capture the randomness. The Krigin model is composed of a mean model, $l(\BFx,\BFP)^{\top} \BFbeta$, and a spatial model, $ M(\BFx,\BFP)$. The mean model captures the global trend, where $l(\BFx,\BFP)$ is a $q \times 1$ vector of basic functions and $\BFbeta$ is a $q \times 1$ vector of trend parameters. The spatial model adopts a covariance function: $\Sigma_F((\BFx,\BFP),(\BFx',\BFP'))=\tau^2 R_F((\BFx,\BFP),(\BFx',\BFP');\theta)$ to quantify the spatial covariance between $M(\BFx,\BFP)$ and $M(\BFx',\BFP')$. Here, $R_F$ is the correlation function with length parameter $\theta$ and variance parameter $\tau$. 

If a Gaussian prior is additionally selected for $\BFbeta$, i.e. $\BFbeta\sim N(\BFb,\Omega)$ with appropriately chosen $\BFb$ and $\Omega$, then stochastic process $F(\BFx,\BFP)$ used to model the mean function is
\begin{align*}
	F(\BFx,\BFP) &\sim \textnormal{Kriging} (l(\BFx,\BFP)^{\top} \BFb, l(\BFx,\BFP)^{\top}\Omega l(\BFx,\BFP) +  \tau^2R_F((\BFx',\BFP'),(\BFx',\BFP');\theta)).
\end{align*}
The metamodel is flexible enough to encompass different types of prior information about the mean response \citep{xie2014bayesian}. The global spatial information can be incorporated by choosing the basis functions $l(\BFx,\BFP)$ and the prior over $\BFbeta$. The local spatial correlation can be represented by the covariance function $\tau^2R_F((\BFx',\BFP'),(\BFx',\BFP');\theta)$.


As mentioned above, we can only get random observations $y(\BFx,\BFP)$ of $f(\BFx,\BFP)$. A common practice to reduce the estimation error is to obtain multiple independent observations of $y(\BFx,\BFP)$. Specifically, for each given $\BFx$ and $\BFP$, we generate $r$ independent values of $\BFu_1, ...,\BFu_r$ from $\BFP$. Then, the $j$-th observations $h(\BFx,\BFu_j)$ is denoted as $y_j(\BFx, \BFP)$. Suppose the current design set contains $n$ design--distribution points denoted  by $\{(\BFx_1, \BFP_1),\cdots,(\BFx_n, \BFP_n) \}$. 
The averaged black-box output at point $(\BFx_i, \BFP_i)$ is $\bar{y}(\BFx_i, \BFP_i)=\frac{1}{r}\sum_{j=1}^{r}y_j(\BFx_i, \BFP_i)$, and  $\bar{\BFY}_n= \{\bar{y}(\BFx_1, \BFP_1), \cdots, \bar{y}(\BFx_n, \BFP_n) \}$. 
If  $\tau$ and $\BFtheta$ are known,  the posterior distribution $F_n(\BFx,\BFP) \triangleq F|\bar{\BFY}_n \sim \textnormal{Kriging}(m_n,k_n)$ can be derived as follows:
\begin{equation*}
		m_n(\BFx,\BFP)= l(\BFx,\BFP)^{\top} \hat{\BFbeta} + \tau^2 R_F((\BFx,\BFP),\cdot;
		{\BFtheta})^T[\tau^2\BFR_F({\BFtheta})+\BFSigma_\epsilon]^{-1} (\bar{\BFY}_n - \BFL^{\top}\hat{\BFbeta}), \label{mean}
\end{equation*}
\begin{equation*}
		\begin{aligned}
			k_n((\BFx,\BFP),({\BFx'},{\BFP'})) &=\tau^2R_F((\BFx,\BFP),({\BFx'},{\BFP'});{\BFtheta}) \\
			&\qquad -\tau^4 \BFR_F((\BFx,\BFP),\cdot;
			{\BFtheta})^T[\tau^2\BFR_F({\BFtheta})+\BFSigma_\epsilon]^{-1}\BFR_F(({\BFx'},{\BFP'}),\cdot;
			{\BFtheta})  \\
			& \qquad  + \eta(\BFx,\BFP)^{\top}\left[\mathbf{\Omega}^{-1}+\BFL[\tau^2\BFR_F({\BFtheta})+\BFSigma_\epsilon]^{-1} \BFL^{\top}\right]^{-1} \eta(\BFx',\BFP')
		\end{aligned}
	\label{cov} 
\end{equation*}
where $\BFL$ is a $q \times n$ matrix that collects the $l(\BFx,\BFP)$ vectors for all $n$ observed points, $\hat{\BFbeta}=(\BFOmega^{-1}+\BFL(\tau^2\BFR_F({\BFtheta})+\BFSigma_\epsilon)^{-1} \BFL^{\top})^{-1}   (L(\tau^2\BFR_F({\BFtheta})+\BFSigma_\epsilon)^{-1} \bar{\BFY}_n+\BFOmega^{-1}\BFb)$ and $\eta(\BFx,\BFP)=l(\BFx,\BFP)-\BFL(\tau^2\BFR_F({\BFtheta})+\BFSigma_\epsilon)^{-1} \tau^2 \BFR_F((\BFx,\BFP),\cdot;
{\BFtheta})$ \citep{rasmussen2004gaussian}.  {The posterior mean $m_n$ approximates $F(\BFx,\BFP)$.  $k_n((\BFx,\BFP),({\BFx'},{\BFP'}))$ is the posterior covariance. 
	In addition, $\tau^2 \BFR_F({(\BFx,\BFP)},\cdot;{})^T$ $=[\text{Cov}[F(\BFx,\BFP),F(\BFx_1,P_1)], \cdots, \text{Cov}[F(\BFx,\BFP),F(\BFx_n, \BFP_n)]]^T$  $\in R^{n\times 1}$, and $\tau^2\BFR_F({\BFtheta})$ $\in R^{n\times n}$ is the covariance matrix of observed  points, and ${\BFSigma_\epsilon} = \text{Diag}\{\frac{1}{r}\sigma_\epsilon^2\,\cdots,\frac{1}{r}\sigma_\epsilon^2\} =\frac{1}{r} \sigma_\epsilon^2 {\BFI}$ $\in R^{n\times n}$.} { If there is no prior knowledge about the trend, which is often the case (as with our experiments in Section \ref{sec3:numericalexperiments}), then a constant value $l(\BFx,\BFP)^{\top} \BFbeta = \beta_0$ is commonly used \citep{xie2014bayesian}. In addition, a vague prior is imposed on $\beta_0$, i.e ${\beta_0} \sim N(0,\Omega)$, where the $\Omega^{-1} \to 0$}. In this case,  the metamodel becomes the nugget-effect model \citep{cressie1992statistics}.

The parameters $\tau$, ${\BFtheta}$ and $\sigma_\epsilon^2$ in the metamodel are unknown and need to be estimated. One common approach is to estimate  $\tau$, and ${\BFtheta}$ with the maximum log-likelihood method \citep{rasmussen2004gaussian} based on the commonly used L-BGFS-B algorithm, and to use the pooled variance of the simulation outputs to estimate $\sigma_\epsilon^2$. Another approach is to impose prior distributions on the parameters and model the parameter uncertainty with a hierarchical Bayesian framework, see \cite{ng2012bayesian}, but this would increase the computational efforts.

$R_F((\BFx,\BFP),(\BFx',\BFP'))$ in the Kriging model needs to be properly chosen to ensure a valid Kriging model. We use the following factorization form {$R_F((\BFx,\BFP),(\BFx',\BFP'))= r_{\mathcal{X}}(\BFx,\BFx') \times r_{\mathcal{P}}(\BFP,\BFP')$ where $r_{\mathcal{X}}$ and $r_{\mathcal{P}}$ are correlation kernels.} We choose the most popular squared exponential correlation kernel for $r_{\mathcal{X}}(\BFx, \BFx^{\prime})=\exp \{-\sum_{i=1}^d
\frac{(\mathrm{x}_i-\mathrm{x}_i^{\prime})^{2}}{2\theta_{1,i}^2}\}$, where $\theta_{1,i}, i = 1, \cdots, d$ are length scale parameters \citep{rasmussen2004gaussian}. Matern correlation kernels and others can also be used. Also, we adopt the form $r_{\mathcal{P}}(\BFP, \BFP^{\prime})=\exp \{- \sum_{j=1}^{l}\frac{D^{2}(P^{(j)}, P^{(j)\prime} )}{2\theta_{2,j}^2}\}$, where $D^{2}(P^{(j)}, P^{(j)\prime} )$ is some measure of the distance between the distributions $P^{(j)}$ and $P^{(j)\prime}$, and $\theta_{2,i}, i = 1, \cdots, l$ are length scale parameters.
In order to make the Kriging model valid, the measure of the distance between distributions $D^2(\cdot,\cdot)$ needs to result in positive definite correlation kernel.
In this paper, we choose { the quadratic Wasserstein distance \citep{bachoc2017gaussian}
	which considers both the mass difference pointwise  and the difference  between points.} Formally, the quadratic Wasserstein distance between two univariate distributions $Q$ and $Q'$, denoted by $WD(Q,Q')$, { is defined by $D^{2}\left(Q, Q^{\prime}\right) = WD(Q,Q') = \textstyle \inf_{\tau\in\gamma(Q,Q')}  E_{(z,z') \sim \tau} [|z-z'|^2$, where $\gamma(Q,Q')$ is set of the joint distributions of $(z,z')$ whose marginal distributions are $Q$ and $Q'$.
	The Wasserstein distance between two univariate distributions can be calculated efficiently with the algorithm provided by \cite{peyre2019computational}.}


\subsection{Approximate model for $g(\mathbf{x})$} \label{sec:gkriging}
Recall that  $F_n(\BFx,\BFP)$ is an approximate model to $f(\BFx,\BFP)$ and aim to derive the model for $g(\BFx)=\text{E}_{\BFP \sim\BFpi}[f(\BFx,\BFP)]$. Consider
\begin{equation*}
	G_n(\BFx)= \text{E}_{\BFP \sim \BFpi}[F_n(\BFx,\BFP)].       \label{G_n}
\end{equation*}
$G_n(\BFx)$ is still a Kriging \citep{DeOliveira2015} with the  mean and covariance functions as follows:
{
	\begin{align}
		\text{E}[G_n(\BFx)]  &= \textstyle \int_{{P}}  \text{E}[ F_n(\BFx, \BFP )] \cdot \pi(\BFP|{\BFD_s})d\BFP = \textstyle \int_{\BFP}   m_n(\BFx, \BFP ) \cdot \BFpi(\BFP|{\BFD_s})d\BFP  \label{equ:G_expectation}
	\end{align}
	\begin{align}
		\text{Cov}[{G}_n(\BFx),{G}_n({\BFx'})] &= \textstyle\int_{\BFP}  \int_{{\BFP'}}  \BFpi(\BFP|{\BFD_s}) \BFpi({\BFP'}|{\BFD_s})  k_n((\BFx, \BFP ),({\BFx'},{\BFP'})) d {\BFP'} d \BFP
		\label{equ:G_cov}
	\end{align}
}
The integrals in (\ref{equ:G_expectation}) and (\ref{equ:G_cov}) can be estimated numerically through Monte Carlo (MC) techniques. Denote $N_{MC}$ as the size of the Monte Carlo samples.  Specifically, we generate $\{\BFP_1,\cdots,{\BFP_{N_{MC}}}\}$ from $\BFpi(\BFP|{\BFD_s})$ to compute
{
	\begin{equation*}
		\text{E}[G_n(\BFx)] \approx \mu_n({\BFx})= \textstyle \frac{1}{N_{MC}}  \sum_{i=1}^{N_{MC}}m_n(\BFx,\BFP_i) \label{MM}
	\end{equation*} 
	\begin{equation*}
			\begin{aligned}
				\text{Cov}[G_n(\BFx),G_n({\BFx'})] &\approx c_n(\BFx,{\BFx'}) = \textstyle\frac{1}{{N_{MC}^2}} \textstyle \sum_{i=1}^{N_{MC} }\textstyle \sum_{j=1}^{N_{MC}}k_n((\BFx,{\BFP_i}),({\BFx'},{\BFP_j})) 
			\end{aligned}
		\label{cc}
	\end{equation*}
	$\text{Var}[G_n(\BFx)] \approx \sigma_n^2(\BFx) = c_n(\BFx,\BFx)$. In fact, we have $\mu_n({\BFx}) = \text{E}\left[\frac{1}{N_{MC}}\sum_{i=1}^{N_{MC}}F_n(\BFx,{\BFP_i})\right]$ and \\ $c_n(\BFx,{\BFx'})= \text{Cov}(\frac{1}{N_{MC}}\sum_{i=1}^{N_{MC} }F_n(\BFx,{\BFP_i}),\frac{1}{N_{MC}}\sum_{j=1}^{N_{MC} }F_n({\BFx'},{\BFP_j}))$. 
	Furthermore, as
	\begin{equation*}
		\hat{G}_n(\BFx)= \textstyle \frac{1}{N_{MC}} \textstyle \sum_{i=1}^{N_{MC} }F_n(\BFx,\BFP_i), \label{MCA}
	\end{equation*}
}
$\hat{G}_n(\BFx)$ is a Kriging with mean  $\mu_n({\BFx})$ and covariance $c_n(\BFx,{\BFx'})$. It can then be used to approximate $G_n(\BFx)$ and used as the approximate model for $g(\BFx)$. 
Popular sampling methods, such as the stick-breaking algorithm, can be used to generate the samples $\BFP_i$ efficiently \citep{hjort2010bayesian}.  The specific procedure for the stick-breaking method \citep{hjort2010bayesian} to sample from the $DP(\alpha, P)$ is as follows:
\begin{enumerate}
	\item	Sample a sequence of $\beta_1, \beta_2,\beta_3, \cdots$ independently from $Beta(1,\alpha)$;
	\item	Sample a sequence of $\theta_1,\theta_2,\theta_3,  \cdots$ independently from the base distribution $P$;
	\item	Set $\alpha_k = \beta_k \Pi_{l=1}^{k-1}(1-\beta_l)$;
	\item	The discrete distribution $\sum_{k=1}^{\infty}\alpha_k\delta_{\theta_k}$ is a sample from the DP.
\end{enumerate}
In practice, we use a truncated discrete distribution $\sum_{k=1}^{K }\pi_k\delta_{\theta_k}$ as the approximated sample from the DP, where $K$ is the minimal integer that satisfies $\sum_{k=1}^{K }\alpha_k >0.9999$ (see \cite{muliere1998approximating} for more details). 

\section{DaBNO-K algorithm}
\label{EGO}
In this section, we propose the DaBNO-K algorithm to solve the DaBNO problem based on the Kriging surrogate model in Section \ref{sec3:algo}. This algorithm uses the idea in the EGO algorithm to guide the subsequent selection of $(\BFx_{n+1},\BFP_{n+1})$. EGO was initially proposed for black-box optimization problems with deterministic outputs. The subsequent evaluation point is selected by maximizing the Expectation Improvement (EI) criterion which takes into account both the prediction of the objective function and the uncertainty associated with such prediction. It can keep a balance of exploitation around the current optimum and exploration in areas with few observed points.

The DaBNO-K algorithm is a modified EGO algorithm. Similar approaches adopting the Knowledge Gradient algorithm \citep{frazier2009knowledge} for the parametric setting was proposed in \cite{Pearce2018} and \cite{toscano2018bayesian}. Two key values are needed when the EI criterion is used. The first one is  the current best value of  $g(\BFx)$ upon which we aim to improve, and the second one is  the predictive distribution of $g(\BFx)$ with a new hypothetical evaluation. The current best value, denoted by $T$, is computed as $\min\{\mu_n({\BFx_1}),\cdots,\mu_n(\BFx_n)\}$, where $\mu_n(\BFx)$ is the predictive mean of the derived Kriging model $\hat{G}_n(\BFx)$, see Section \ref{sec:gkriging}. Denote the predictive distribution $\hat{G}_{n+1}(\BFx)$ conditional on the evaluation at an arbitrary point $(\BFx_{n+1},\BFP_{n+1})$ as $\hat{G}_{n+1}(\BFx|{\BFx_{n+1}},\BFP_{n+1})$, then
\begin{equation}
		\textstyle \hat{G}_{n+1}(\BFx|{\BFx_{n+1}},\BFP_{n+1})=\frac{1}{N_{mc}}\sum_{i=1}^{N_{mc}}F_{n+1}((\BFx, \BFP_i)|\BFx_{n+1},\BFP_{n+1}).
	\label{sum}
\end{equation} 
where $F_{n+1}((\BFx, \BFP_i)|\BFx_{n+1},\BFP_{n+1})$ is the  predictive distribution  of $F_{n+1}(\BFx,\BFP)$ conditional on evaluation at $(\BFx_{n+1},\BFP_{n+1})$.
Denote the value evaluated at $(\BFx_{n+1},\BFP_{n+1})$ by $y_{n+1}$.  While $(\BFx_{n+1},\BFP_{n+1})$ are fixed, $y_{n+1}$ are  normally distributed random variables.
We can then derive
$F_{n+1}((\BFx,\BFP)|\BFx_{n+1},\BFP_{n+1},y_{n+1}) \sim \textnormal{Kriging}(m'_{n+1},k'_{n+1})$ as follows:
\begin{equation*}
		\begin{aligned}
			m'_{n+1} (\BFx, \BFP) &=  \textstyle m_n(\BFx,\BFP) + \frac{k_n((\BFx,\BFP),(\BFx_{n+1},\BFP_{n+1}))}{k_n((\BFx_{n+1},\BFP_{n+1}),(\BFx_{n+1},\BFP_{n+1})) +\frac{1}{r}\sigma_\epsilon^2}[y_{n+1}-m_n(\BFx,\BFP)] \\
			& =  \textstyle m_n(\BFx,\BFP) + \frac{k_n((\BFx,\BFP),(\BFx_{n+1},\BFP_{n+1}))}{\sqrt{k_n((\BFx_{n+1},\BFP_{n+1}),(\BFx_{n+1},\BFP_{n+1})) +\frac{1}{r}\sigma_\epsilon^2}}Z,
		\end{aligned}
\end{equation*}
\begin{equation}
		\begin{aligned}
			k'_{n+1}((\BFx,\BFP),(\BFx',\BFP')) &= \textstyle k_n((\BFx,\BFP),(\BFx',\BFP'))   -   \frac{ k_n((\BFx,\BFP),(\BFx_{n+1},\BFP_{n+1}))k_n((\BFx',\BFP'),(\BFx_{n+1},\BFP_{n+1}))}{k_n((\BFx_{n+1},\BFP_{n+1}),(\BFx_{n+1},\BFP_{n+1})) +\frac{1}{r}\sigma_\epsilon^2},
		\end{aligned}
\end{equation}
where $Z \sim {N}(0,1)$.

Then from Equation (\ref{sum}), we can  derive $\hat{G}_{n+1}(\BFx|{\BFx_{n+1}},\BFP_{n+1},y_{n+1})$ with  mean and covariance as follows:
\begin{equation*}
	\resizebox{\hsize}{!}{$
		\begin{aligned}
			\mathrm{E}\left[\hat{G}_{n+1}(\BFx|{\BFx_{n+1}},\BFP_{n+1},y_{n+1})\right] 
			&= \textstyle \frac{1}{N_{mc}} \sum_{i=1}^{N_{mc}} m'_{n+1} (\BFx, \BFP_i)  \\
			&= \textstyle \frac{1}{N_{mc}} \sum_{i=1}^{N_{mc}} m_n(\BFx, \BFP_i)  + Z\frac{1}{N_{mc}} \sum_{i=1}^{N_{mc}}  \frac{k_n((\BFx, \BFP_i),(\BFx_{n+1},\BFP_{n+1}))}{\sqrt{k_n((\BFx_{n+1},\BFP_{n+1}),(\BFx_{n+1},\BFP_{n+1})) +\frac{1}{r}\sigma_\epsilon^2}}\\
			&\sim \textstyle N\left(\frac{1}{N_{mc}} \sum_{i=1}^{N_{mc}} m_n(\BFx, \BFP_i) , \left(\frac{1}{N_{mc}} \sum_{i=1}^{N_{mc}}  \frac{k_n((\BFx, \BFP_i),(\BFx_{n+1},\BFP_{n+1}))}{\sqrt{k_n((\BFx_{n+1},\BFP_{n+1}),(\BFx_{n+1},\BFP_{n+1})) +\frac{1}{r}\sigma_\epsilon^2}}\right)^2\right),
		\end{aligned}
		$}
\end{equation*}
\begin{equation*}
	\resizebox{1\hsize}{!}{$
		\begin{aligned}
			&\mathrm{Cov} [\hat{G}_{n+1}(\BFx|{\BFx_{n+1}},\BFP_{n+1},y_{n+1}), \hat{G}_{n+1}(\BFx'|{\BFx_{n+1}},\BFP_{n+1},y_{n+1})] = \textstyle \frac{1}{N_{mc}^2} \sum_{i=1}^{N_{mc}}  \sum_{j=1}^{N_{mc}} k'_{n+1}((\BFx, \BFP_i),(\BFx',\BFP_j)).
		\end{aligned}
		$}
\end{equation*}


Then we have $\hat{G}_{n+1}(\BFx|{\BFx_{n+1}},\BFP_{n+1}) \sim N\left(\mu_n'(\BFx), \sigma_n'^2(\BFx)\right)$, where 
\begin{equation*}
	\mu_n'(\BFx) = \textstyle \frac{1}{N_{mc}} \sum_{i=1}^{N_{mc}} m_n(\BFx, \BFP_i) ,
\end{equation*}
\begin{equation*}
	\resizebox{1\hsize}{!}{$
		\begin{aligned}
			\sigma_n'^2(\BFx) &= \textstyle \left(\frac{1}{N_{mc}} \sum_{i=1}^{N_{mc}}  \frac{k_n((\BFx, \BFP_i),(\BFx_{n+1},\BFP_{n+1}))}{\sqrt{k_n((\BFx_{n+1},\BFP_{n+1}),(\BFx_{n+1},\BFP_{n+1}))
					+\frac{1}{r}\sigma_\epsilon^2}}\right)^2  + \frac{1}{N_{mc}^2} \sum_{i=1}^{N_{mc}}  \sum_{j=1}^{N_{mc}} k'_{n+1}((\BFx, \BFP_i),(\BFx, \BFP_j)).
		\end{aligned}
		$}
\end{equation*}

The proposed EI  criterion is
\begin{equation}
		\begin{aligned}
			\text{EI}_T({\BFx_{n+1}},\BFP_{n+1})&= \text{E}_{\hat{G}_{n+1}}[(T-\hat{G}_{n+1}(\BFx_{n+1}|{\BFx_{n+1}},\BFP_{n+1}))^+|\bar{\BFY}_n] 
			\\
			& =\Delta \Phi\left(\frac{\Delta}{\sigma_{n}^{\prime}\left(\BFx_{n+1}\right)}\right)+\sigma'_{n}\left(\BFx_{n+1}\right) \phi\left(\frac{\Delta}{
				\sigma_{n}^{\prime}\left(\BFx_{n+1}\right)}\right)
		\end{aligned}
	\label{equ:EGOIU}
\end{equation}
where $\Delta = T-\mu_{n}^{\prime}\left(\BFx_{n+1}\right)$. Although the value of $g(\BFx_{n+1})$ is unknown, we have derived an approximation given by the distribution $\hat{G}_{n+1}(\BFx_{n+1}|\BFx_{n+1},\BFP_{n+1})$. Hence, the improvement given by $T-\hat{G}_{n+1}(\BFx_{n+1}|{\BFx_{n+1}},\BFP_{n+1})$ is averaged with respect to  $\hat{G}_{n+1}(\BFx_{n+1}|\BFx_{n+1},\BFP_{n+1})$. 

The DaBNO-K algorithm starts with initializing the surrogate model with ${n_0}$ initial points (See Initialization step in Algorithm 1). In each of the following iterations, the next point to evaluate is selected as $\arg\max_{(\BFx_{n+1},\BFP_{n+1})} \text{EI}_T(\BFx_{n+1},\BFP_{n+1})$. We then conduct
$r$ replications of function evaluations to obtain $\bar{y}(\BFx_{n+1},\BFP_{n+1})$ and update the Kriging model. This process is iterated until it reaches the computing budget $N$ (the total number of function evaluations). Denote $\mathcal{V} = \{\BFx_1,\cdots, \BFx_n\}$ as the set of evaluated decision vectors, then $\hat{\BFx}^*_N=\arg\min_{\BFx\in \mathcal{V}} \mu_n(\BFx)$ is  the final estimated minimizer.  Algorithm \ref{algo} provides an outline of the proposed algorithm.


\begin{algorithm}
	\caption{the DaBNO-K algorithm}
	\label{algo}
	\begin{algorithmic}[1]
		\State Given real world data $\BFD_s= \{\BFu_1,\cdots,\BFu_s\}$  and a  prior $DP(\alpha, \BFP^0)$ on the distribution $\BFP$, derive the posterior $\BFpi(\BFP|{\BFD_s})$;
		\State \textbf{Initialization}:
		Generate $ \{{\BFx_1},\cdots,\BFx_{n_0}\}$ using Latin Hypercube Sampling with support $\mathcal{X}$ and generate  $\{\BFP_1,\cdots,\BFP_{n_0}\}$   randomly from $\BFpi(\BFP|{\BFD_s})$ to obtain an initial sample points: $\{({\BFx_1},\BFP_1),\cdots,(\BFx_{n_0},\BFP_{n_0})\}$. Run $r$ function evaluations at each point to obtain $\bar{\BFY}_{n_0}=[\bar{y}({\BFx_1},\BFP_1),\cdots,\bar{y}(\BFx_{n_0},\BFP_{n_0})]^T$.
		Set $\mathcal{V} = \{\BFx_1,\cdots, \BFx_{n_0}\}$;
		\State \textbf{Surrogate model}: Construct a  Kriging model $F_{n_0}(\BFx,\BFP)$ based on $\bar{\BFY}_{n_0}$; 
		\While{$ nr \leq N -1 $} 
		\State Generate $\BFP_i, i= 1, \cdots,{{N_{MC}}}$ from $\BFpi(\BFP|{\BFD_s})$; Derive  $\hat{G}_n(\BFx,\BFP)\sim \text{Kriging}(\mu_{n},c_{n})$ for $g(\BFx)$;
		\State  \textbf{Selection:} Choose the subsequent  point $({\BFx_{n+1}},{\BFP_{n+1}})$ that maximizes $ \text{EI}_T(\BFx_{n+1},\BFP_{n+1})$ in Equation \eqref{equ:EGOIU};
		\State Run $r$ evaluations at $({{\BFx_{n+1}}},{\BFP_{n+1}})$ to obtain  $\bar{y}({\BFx_{n+1}},{\BFP_{n+1}})$, set ${\bar{\BFY}_{n+1}} = [{\bar{\BFY}_{n}}, \bar{y}({\BFx_{n+1}},{\BFP_{n+1}})]^T$;
		\State \textbf{Update}: update the  Kriging model $F_{n+1}(\BFx,\BFP)$ given ${\bar{\BFY}_{n+1}}$;
		\State Set $\mathcal{V} = \mathcal{V} \cup \{\BFx_{n+1}$\} and $n = n+1$;
		\EndWhile
		\State \textbf{return} $\hat{\BFx}^*_N=\arg\min_{\BFx\in \mathcal{V}} \mu_n(\BFx)$.
	\end{algorithmic}
\end{algorithm}

\section{Convergence analysis of the DaBNO-K algorithm} 
\label{convergence}
In this section, we aim to study the convergence of the DaBNO-K algorithm. Here,  we aim to show that 1) the objective function in Equation \eqref{objfunction} can be well-solved by our proposed algorithm when the amount of real world data points $s$ is fixed and the computing budget $N$ goes to infinity;  2) if both $s$ and $N$ goes to infinity, the estimated optimal value/solution would converge to the optimal value/solution of the true problem in Equation \eqref{trueobjective}.
We first introduce our main assumptions.

\begin{assumption}\label{assum2}
	$\mathcal{X}$ is a compact space.
\end{assumption}
\begin{assumption}\label{assum3}
	The parameters $\tau, \bm{\theta}, \sigma_\epsilon$ of the Kriging model are assumed known.
\end{assumption}
Recall  $\BFx_s =  \arg\min g(\BFx)$ and let $\widehat{\BFx}_N^*$ be the solution returned by the proposed algorithm when the simulation budget is $N$. The main theorem we are going to prove is as follows:
\begin{theorem}\label{thm:main theorem}
	Suppose Assumptions \ref{assum2}--\ref{assum3} hold.  Then as $N$ tends to infinity, we have
	\begin{equation*}
		g(\widehat{\BFx}_N^*)-g(\BFx_s)\leq o_p(1),
	\end{equation*}
	i.e., $\widehat{\BFx}_N^*$ nearly minimizes $g(\BFx)$.
\end{theorem}

The proof of this theorem is decomposed into two main steps. In the first step, we aim to show that the posterior mean converges to the true objective:
\begin{equation}\label{eq: mu_n to g}
	\left|\mu_n(\BFx)-g(\BFx)\right|\to_p 0
\end{equation}
uniformly for all $\BFx$.
In the second step, we are going to show that the points we visit in our algorithm are dense, and our idea follows from \cite{locatelli1997bayesian}.
These two results then naturally imply the convergence of our algorithm.
A direct corollary of the above theorem is that the sequence of points returned by Algorithm \ref{algo}  converges to the global minimizer in probability if the set of optimum points includes only one element.

\begin{corollary}\label{lem: consistency of xN to xg}
	Suppose Assumptions \ref{assum2}--\ref{assum3} hold.  If the global optimal solution of $g(\BFx)$ is unique, then as $N \to \infty$,
	\begin{equation*}
		\hat{\BFx}_N^*\to_p \BFx_s.
	\end{equation*}
\end{corollary}
When both $s$ and $N$ go to infinity, we have the following result.
\begin{corollary}\label{corollary:2}
	Suppose conditions in Corollary \ref{lem: consistency of xN to xg} hold.   Then for all $\epsilon>0$, there exists $s_0 >0$  such that for all $s>s_0$, there exists $N_0$ which depends on $s$, such that  $\mathbb{P}(|\hat{\BFx}_N^*-\BFx^*|\geq \epsilon)<\epsilon$ and  $P(| g(\hat{\BFx}_{N}^*)-f(\BFx^*,\BFP^c)|\geq \epsilon)<\epsilon$ 
	for any $N\geq N_0$.
	In other words, as $s$ and $N$ tends to infinity,
	\begin{equation*}
		\hat{\BFx}_N^* \to_p \BFx^*, \quad g(\hat{\BFx}_{N}^*) \to_p f(\BFx^*,\BFP^c).
	\end{equation*}
\end{corollary}

\begin{remark}
	Corollary \ref{corollary:2} shows that, when both  $s$ and  $N$ become large enough, the optimal solution/value returned from our algorithm can converge to the optimal solution/value of the true problem.
\end{remark}

\section{Experiments} \label{sec3:numericalexperiments}
In this section, to assess the convergence properties of the DaBNO objective function and the proposed DaBNO-K algorithm, we first test their empirical convergence with a classic inventory problem \citep{fu1997techniques,jalali2017comparison} in Section \ref{experi1}. In Section \ref{experi2}, we further compare the finite sample performance of DaBNO-K with several classical approaches. While the primary focus of the paper's method lies in continuous optimization, it is important to highlight that our approach extends its applicability to spaces accommodating discrete variables. This versatility is evident in our comparison studies, where we address both continuous problems (including two synthetic functions and the inventory problem above) and a discrete critical care facility allocation problem. All experiments in this section are run on a high-performance-computing cluster with the following specifications: CPU E5-2690 v3  @2.60GHz with 24 cores.

\subsection{Empirical convergence results} \label{experi1}
In this section, we numerically test the convergence results in this work. These include the convergence of DaBNO objective function (in Section \ref{sec3:properties}) and the convergence of DaBNO-K algorithm (in Section \ref{convergence}). Recall that $s$ is the number of real-world observations of $\BFu$ and $N$ is the computing budget (number of function evaluations). More specifically, our numerical results are
\begin{itemize}
	\item the empirical convergence of the DaBNO objective function to the true objective function (when $s\to \infty$, Section \ref{e1});
	\item the empirical convergence of the proposed DaBNO-K algorithm in solving the DaBNO problem (when $s$ fixed and $N \to \infty$, Section \ref{e2});
	\item the empirical convergence of the proposed DaBNO-K algorithm in solving the true objective function (when $s$ fixed and $N \to \infty$, Section \ref{e3}).
\end{itemize}
These results are obtained with a widely used \textbf{inventory problem} \citep{fu1997techniques,jalali2017comparison}.

\subsubsection{Setup of the inventory problem}
We consider a continuous inventory problem where a company periodically adjusts its inventory with the following rules  \citep{fu1997techniques}. There are two key values that need to be specified in the inventory problem: the basic ordering level denoted by $\mathrm{x}_1$ and the order-up-to level $\mathrm{x}_2$.  The levels $\mathrm{x}_1$, $\mathrm{x}_2$ and the demand level are all continuous values. Such setting is applicable for scenarios when the good are uncountable like the oil in a gas station. Denote the inventory position as $p$.
At the beginning of each period, for example, each week, if $p <\mathrm{x}_1$, the company will order $\mathrm{x}_2-p$ items. If  $p \geq \mathrm{x}_1$, no action is needed. The decision variable in the inventory problem is $\BFx=[\mathrm{x}_1,\mathrm{x}_2]$ and $\{\BFx=[\mathrm{x}_1,\mathrm{x}_2]|\mathrm{x}_1 \in [10000,22500], \mathrm{x}_2 \in [22600,35000]\}$. It is assumed that the customer demand, denoted by $u$, for each period, is i.i.d following a distribution denoted by $P$. Other hyper-parameters including { fixed ordering cost, unit cost, holding cost, and backorder cost} are set to be 100, 1, 1 and 100 respectively \citep{jalali2017comparison}. Each evaluation run has a length of 1000 periods and the first 100 periods are discarded as a warm-up stage. Given  $\BFx$ and $P$,  { $y(\BFx, P) = h(\BFx,u), u \sim P$ is the steady-state cost averaged over the 900 periods and we denote expected output as $f(\BFx, P)$}.

The true distribution of $u$ is assumed to be $P^c = \exp(\lambda^c= 0.0002)$. The expected cost function has an analytical form:
\begin{equation}
		f(\BFx, {P}^c = \textstyle \mbox{exp}(\lambda^c))=\frac{1}{100}\left[\frac{1}{\lambda^c}+\frac{100+\mathrm{x}_1-\frac{1}{\lambda^c}+0.5 \lambda^c({\mathrm{x}_2}^2-{\mathrm{x}_1}^2)+\frac{101}{\lambda^c}e^{-\lambda^c \mathrm{x}_1}}{1+\lambda^c(\mathrm{x}_2-\mathrm{x}_1)} \right]. 
	\label{analytical}
\end{equation}	
However,  $P^c$ is unknown in practice and only $\bm{D}_{s}=\{u_1, \cdots, u_{s}\}$, where $u_i \sim P^c,  i = 1, \cdots,s$ can be observed. Therefore, the analytical expression in (\ref{analytical}) can not be used to calculate the optimal decision.  We can only infer the black-box surface of (\ref{analytical}) through the simulator from which $y(\BFx, P)$ can be observed for a given ($\BFx$, $P$) pair.

Adopting the DaBNO framework, we model $P^c$ as a random variable $P$. {We use  a  weakly informative Dirichlet process prior  with $\alpha = 1$ and $P_0 = \text{uniform} (0, \max({u_1, ...u_s}))$ (as recommended in \cite{gelman2013bayesian}).
	We aim to solve  $\min_{\BFx} g(\BFx) = \text{E}_{P \sim \pi} [{f}(\BFx,P)]$} based on our proposed Kriging-based algorithm.


\subsubsection{Empirical Convergence and Asymptotic Normality of the DaBNO} \label{e1}
\begin{figure}[!t]
	\begin{minipage}{0.49\linewidth}
		\centering
		\includegraphics[width=1\linewidth]{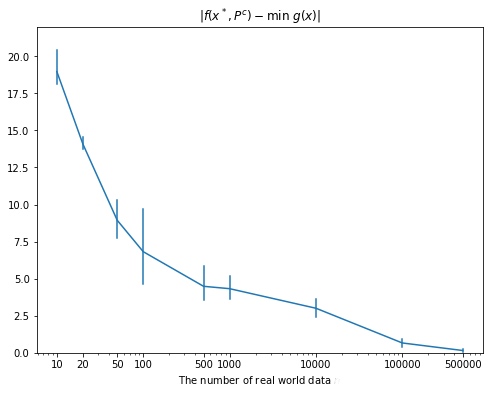}
	\end{minipage}
	\begin{minipage}{0.49\linewidth}
		\centering
		\includegraphics[width=1\linewidth]{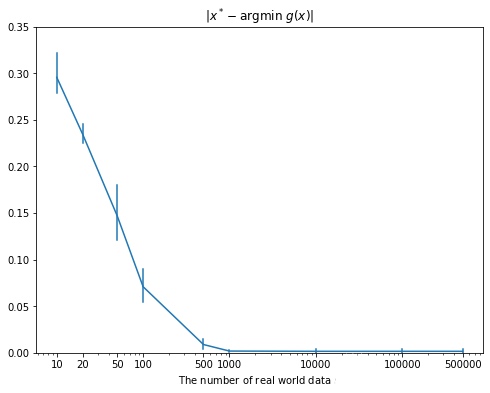}
	\end{minipage}
	\caption{Empirical convergence of the optimal value/solution of the DaBNO objective to the true optimal value/solution as $s$ increases}
	\label{fig:fconvergence}
\end{figure}

We aim to illustrate the empirical convergence of the optimal value/solution of the DaBNO to the optimal value/solution of the $f(\BFx, P^c)$ for finite $s$. We consider two evaluation metrics,  $|f(\BFx^*,P^c)- \min g(\BFx)|$ and $||\BFx^* - \arg\min g(\BFx)||$.  We approximate the unknown  $g(\BFx)$ by {$\frac{1}{10000} \sum_{i=1}^{10000} \bar{y} (\BFx, P_i)$,  where $\bar{y}(\BFx, P_i)$ is the sample mean output} over 1000 replications at $(\BFx, P_i)$.  We then obtain $\min_x g(\BFx)$ and $\arg\min_xg(\BFx)$ with {a grid search over a large set of the design points}. We set  $s$ to be 10, 20, 50, 100, 500,  1000, 10000, 100000 and 500000. We conducted 100 trials and calculated the mean values of the $|f(\BFx^*,P^c)- \min g(\BFx)|$ and $||\BFx^* - \arg\min g(\BFx)||$  and their 95\% confidence intervals over the 100 replications. 
The results are displayed in Figure \ref{fig:fconvergence} which shows that  the optimal value/solution of the DaBNO  approach  the optimal value/solution of  $f(\BFx, P^c)$ as $s$ becomes large, validating Theorem \ref{theorem:1}  $(ii)$ and $(iii)$ for finite $s$.

In addition, the Jarque–Bera normality test \citep{jarque1980efficient} on the values of $\sqrt{s} (f(\BFx^*,P^c)- \min g(\BFx))$ and Mardia’s multivariate normality test  \citep{mardia19809} on the values of $\sqrt{s}(\BFx^* - \arg\min g(\BFx))$ over the 100 replications  both confirm normality for $s = 500000$, which helps validate Theorem \ref{theorem:2} $(ii)$ and $(iii)$ for finite $s$.

\subsubsection{Empirical Convergence of the Proposed Algorithm to the DaBNO problem} \label{e2}
In this section, we examine the empirical convergence of our proposed algorithm in solving the DaBNO problem (Theorem \ref{thm:main theorem} and Corollary \ref{lem: consistency of xN to xg}).
We define gGAP = $|\min_{\BFx} g(\BFx) - g(\hat{\BFx}^*)|$ and gxGAP = $||\arg \min_{\BFx} g(\BFx) - \hat{\BFx}^*||$  to test the empirical convergence of the returned value/solution of the proposed DaBNO-K algorithm to the optimal value/solution of $g(\BFx)$.  

\begin{figure}[!t]
	\centering
	\includegraphics[width=1\linewidth]{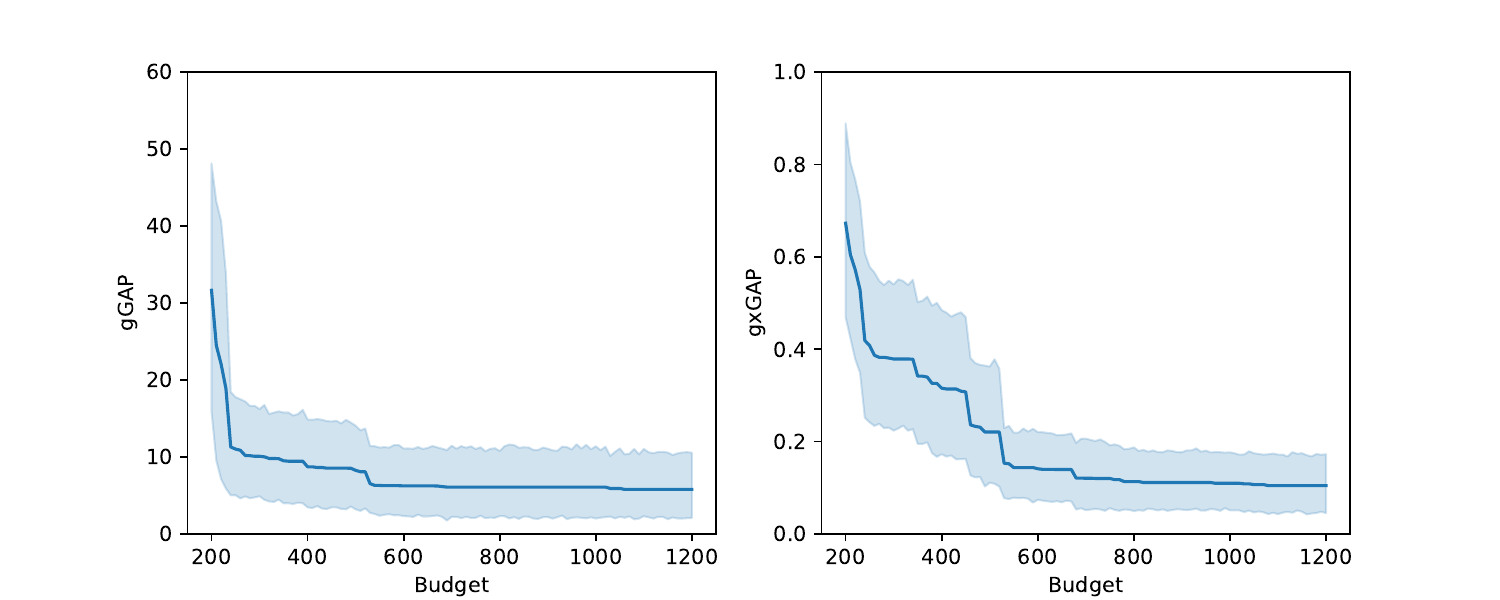}
	\caption{ Empirical performance in terms of gGAP and gxGAP of our surrogate based algorithm as the size of the budget increases. The number of real-world data $s$ is $1000$.}
	\label{fig:gconvergence}
\end{figure}

The initial size of the design sample is chosen to be 20 ($10\times$ dimension of the function), as recommended in \cite{jones1998efficient}. An additional 100 points are sequentially selected according to the proposed algorithm.
$r = 10$ replications are conducted at each selected point.  
The total budget is then $(20+100) \times 10 = 1200$. We conducted 100 trials and calculated the means of both gGAP and gxGAP and their 95\% confidence intervals. We set $s=1000$.
Figure \ref{fig:gconvergence} which shows that as the number of iterations increases, the optimal value/solution returned by our algorithm can get close to the optimal value/solution of the $g(\BFx)$. We could further improve the performance by either increasing the  budget or increasing the number of Monte Carlo approximations $N_{MC}$  (currently $N_{MC} = 100$).

\subsubsection{Empirical Convergence of the Proposed Algorithm to the true objective function} \label{e3}
\begin{figure} [!t]
	\centering
	\includegraphics[width=1\linewidth]{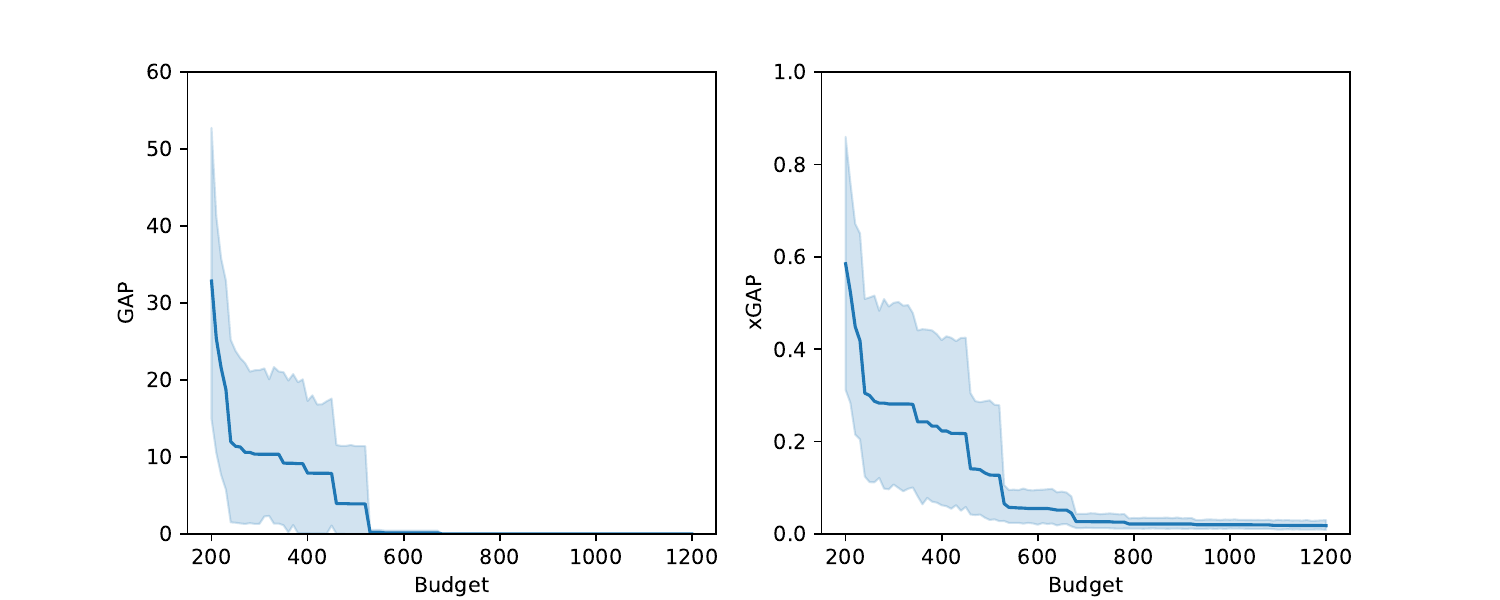}
	\caption{Empirical performance in terms of GAP and xGAP of our surrogate based algorithm as the size of the budget increases. The number of real-world data $s$ is $1000$.}
	\label{fig:convergence}
\end{figure}

In this section, we examine the empirical convergence of our proposed algorithm in solving the true problem.
We consider two evaluation metrics: GAP = $|f(\bm{x^*}, P^c)-{f}(\hat{\BFx}^*,P^c)|$ and xGAP = $||\BFx^*-\hat{\BFx}^*||$. 
The results are displayed in  Figure \ref{fig:convergence}.

It shows that the optimal value and optimal solution returned from our algorithm will empirically converge (with a reasonably finite $s$ and $N$) to the true optimal value and true optimal solution, validating also the convergence properties provided in Corollary \ref{corollary:2} for finite $s$.

\subsection{Comparison of the DaBNO-K with existing approaches}  \label{experi2}
In this section, we aim to compare our DaBNO-K approach with other commonly used approaches. We will first test the results on two synthetic functions and then apply the results on two practical problems --  an inventory example and a critical care facility allocation problem. In all experiments, the true underlying uncertainty distributions are multimodal and complex, making it hard to choose appropriate parametric forms.

\subsubsection{Alternative approaches}
As was mentioned in Section \ref{sec3:introduction}, it is common to ignore the uncertainty of the distribution and solve  (\ref{approproblem}) using an estimated  distribution $\hat{P}$, typically the empirical cumulative distribution. The approach is referred to as the \textbf{hist approach}. Another popular approach is  the \textbf{parametric approach}, which assumes a parametric family for the distribution and estimates the parameters from real world data. The hist and parametric approach then can be easily solved using sequential Kriging optimization (SKO), a variant of the EGO algorithm \citep{Huang2006}. 

\subsubsection{Synthetic problems}
\paragraph{\textbf{Set up of the numerical problems}}
We consider two synthetic functions for comparison: Griewank and StybTang, both with three dimensions. The two functions are commonly used in literature for evaluating optimization algorithms. StybTang has a relatively flat surface, while Griewank's suface are highly volatile. For each function, we choose the first two variables as the design parameter $\bm{x}$ and another one as the uncertain input $\BFu$. The true distribution of $\BFu$ is a mixture of a lognormal distribution with a mean of 10 and a standard deviation of 10, and a lognormal distribution with a mean of 20 and a standard deviation of 5, with a proportion of 0.5 and 0.5. We choose the mixture of lognormal distribution to mimic a complex multi-modal uncertain distribution, which is very common in reality. The analytical expressions of the two functions are shown in the appendix.

For these experiments, we compare four approaches, DaBNO-K, hist, and two parametric methods which assume the uncertain distribution follows an exponential or a lognormal distribution, denoted as parametric-exponential and parametric-lognormal. The four approaches begin with the same twenty initial design points in each trail, and subsequently select forty points. For each design point, the observations are taken as the average of two evaluations ($r=2$). The total budget is 120 ($(20+40) \times 2$). The true optimal values $f(\bm{x}^*, P^c)$ are estimated through extensive monte carlo sampling as -1.85 and -0.02 for the Griewank function and the StybTang function respectively. 

\paragraph{\textbf{Performance}} For each pair of algorithm and test function, we repeat the optimization process for 100 trails under three different values of $s$ (10, 100, 1000). As $s$ is the number of real-world observations of $\BFu$, these three scenarios represent different uncertainty levels for the distribution of $\BFu$. We use GAP = $| f(\BFx^*, P^c)-{f}(\hat{\BFx}^*,P^c)|$ as the evaluation metric. 

\begin{figure}
	\begin{minipage}{0.49\linewidth}
		\centering
		\includegraphics[width=1\linewidth]{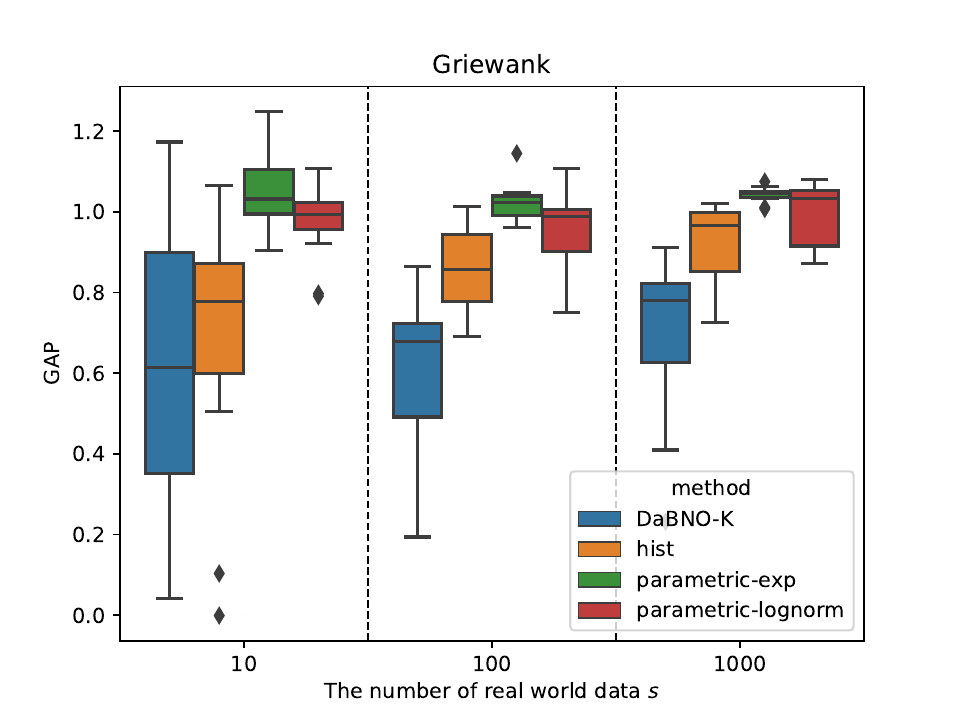}
	\end{minipage}
	\begin{minipage}{0.49\linewidth}
		\centering
		\includegraphics[width=1\linewidth]{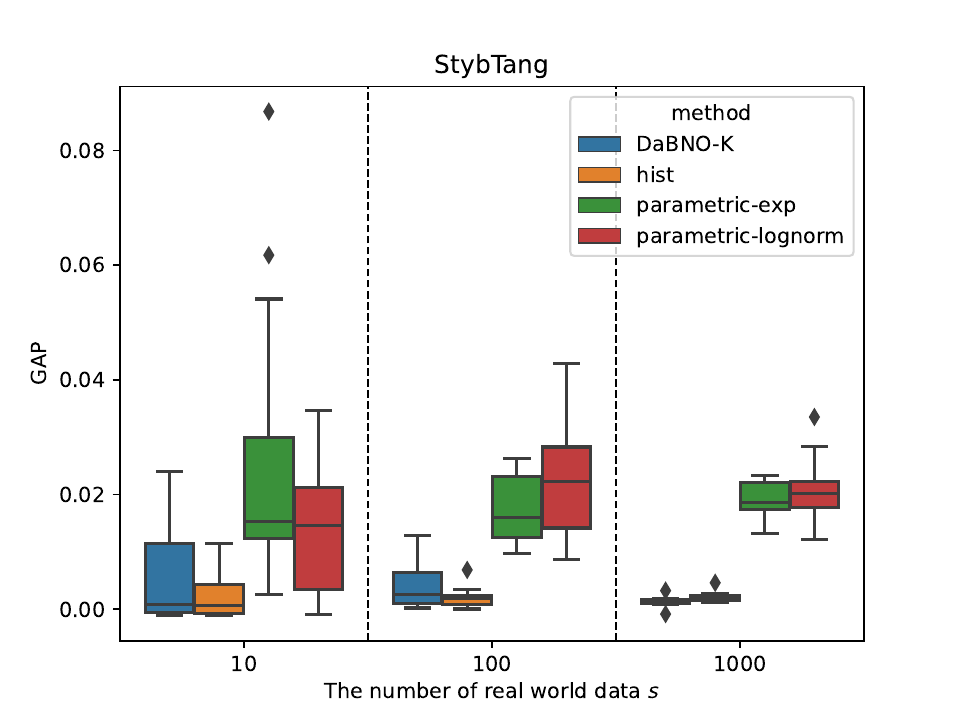}
	\end{minipage}
	\caption{The Griewank  and StybTang function : the boxplot of the GAP values for three approaches (DaBNO-K, hist, parametric-exp and parametric-lognorm) under different uncertainty levels of the distribution of $\BFu$}
	\label{fig:numerical}
\end{figure}

We plot the boxplot of the GAP in the 100 trails in Figure \ref{fig:numerical}. The DaBNO-K approach is more robust to different levels of uncertainty in the distribution of $\BFu$, compared to the hist and parametric approaches across different functions. For the parametric approaches which assume a wrong uncertain distribution (the parametric-exp and parametric-lognormal), the returned optimum deviates from the true optimum, and this deviation cannot be remedied by increasing the number of real world data.

\subsubsection{Practical problems}
\paragraph{\textbf{Set up of the inventory problem}}
We consider the same inventory example as above, but now assume the true uncertain distribution to be a mixture of log-normal distributions instead of a simple exponential distribution. The true distribution is a  mixture of a lognormal distribution with a mean of 5000 and a standard deviation of 5000,  and a lognormal distribution with a mean of 10000 and a  standard deviation of 5000, with a proportion of 0.5 and 0.5. For this experiment, we compare four approaches, DaBNO-K, hist, and two parametric methods which assume the uncertain distribution follows an exponential or a lognormal distribution, denoted as parametric-exponential and parametric-lognormal.  The four approaches begin with the same twenty initial design points for each replication, and subsequently select forty points. For each design point, the observations are taken as the average of ten evaluations. The total budget is 600 ($(20+40) \times 10$).   The estimated true optimal solution $\BFx^* $ is $[22195, 26398]$ and the estimated optimal value $f(\BFx^*, P^c)$  is 320.8. 

\begin{figure}
	\centering
	\includegraphics[width=0.7\linewidth]{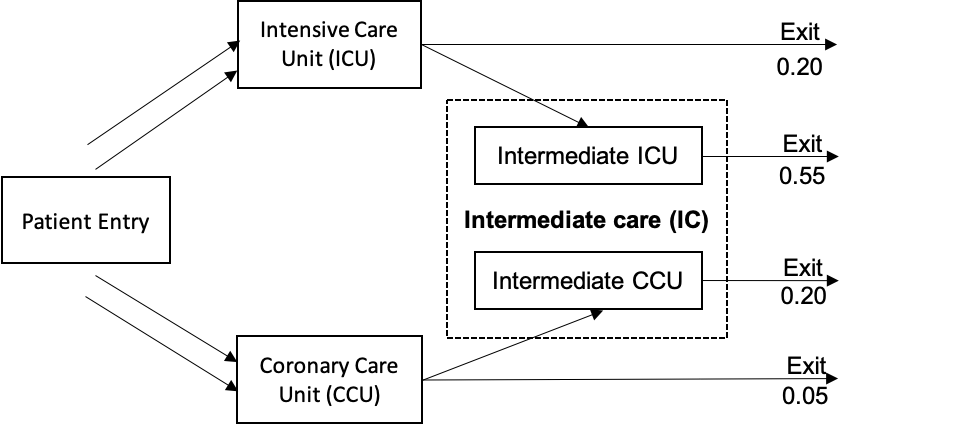}
	\caption{Critical care facility}
	\label{fig:ccf}
\end{figure}

\paragraph{\textbf{Set up of the critical care facility problem}}
In this subsection, we use a critical care facility allocation problem \citep{ng2006reducing, xie2014bayesian} to further illustrate the performance of the proposed algorithm. The output of the system is the steady-state number of patients denied a bed per day. Patients arriving at the facility will be allocated, depending on their health condition, to either the intensive care unit (ICU) or the coronary care unit (CCU), and then exit the facility or go to the intermediate ICU (IICU) or the intermediate CCU (ICCU). The IICU and the ICCU share the beds in the Intermediate Care (IC). Arriving patients who cannot get a bed in ICU or CCU will be denied entry. If one patient is supposed to move to the IC but there is no bed available in the IC, he/she will stay put (in ICU/CCU) and queue for a bed in the IC. When a bed in the IC becomes available, the first patient in the queue will be transferred.  Figure \ref{fig:ccf} depicts the system.

There are three design variables in the system: the number of beds in the ICU, the CCU and the IC, denoted by $\BFx = [\mathrm{x}_1,\mathrm{x}_2,\mathrm{x}_3]^T$. It can be expected that, the larger the number of beds in the three units, the lower the expected number of patients will be denied entry per day. However, the resources in the facility will be always limited. The difficulty is to balance the number of beds of the three different units under limited resources. The cost of holding an IC bed is usually half of that of the ICU and CCU, as fewer staff are needed per patient in the IC \citep{plate2017utilisation}. Therefore, we consider the design space to be $\mathcal{X} = \{\BFx = [\mathrm{x}_1,\mathrm{x}_2,\mathrm{x}_3]| 26 \leq \mathrm{x}_1+\mathrm{x}_2+ 0.5\mathrm{x}_3 \leq 28, \mathrm{x}_1 \in N_+, \mathrm{x}_2 \in N_+, \mathrm{x}_3 \in N_+ \}$. The total number of candidate design points is 1273.

The system also includes six input processes: the inter-arrival time of patients, the stay duration in ICU, CCU, IICU and ICCU, and the routing of patients. The setting for the true uncertain distribution, denoted by $\bm{P^c} = [P^{c(1)},P^{c(2)},P^{c(3)},P^{c(4)},P^{c(5)},P^{c(6)}]$, is as follows. The inter-arrival time of patients ($P^{c(1)}$) is exponentially distributed with the rate being 3.3/day, i.e. patients are arriving according to a Poisson process. The stay duration in ICU, CCU, IICU and ICCU all follow mixtures of lognormal-distributions: 
$P^{c(2)}$ follows a mixture of two lognormal distributions (proportion being 0.5 and 0.5) with a mean of 3.4 and 6.8 and a standard deviation of 3.5;
$P^{c(3)}$ follows a mixture of two lognormal distributions (proportion being 0.5 and 0.5) with a mean of 3.8 and 7.6 and a standard deviation of 1.6;
$P^{c(4)}$ follows a mixture of two lognormal distributions (proportion being 0.5 and 0.5) with a mean of 15 and 30 and a standard deviation of 7;
$P^{c(5)}$ follows a mixture of two lognormal distributions (proportion being 0.5 and 0.5) with a mean of 17 and 34 and a standard deviation of 3.
The routing of each patient $(P^{c(6)})$ follows a discrete distribution with probability 0.2, 0.55, 0.2 and 0.05. 

We consider the following three approaches in the comparison: the DaBNO, the hist, the parametric approach which assume $P^{c(1)}$ to be exponentially distributed, $P^{c(2)}, P^{c(3)}, P^{c(4)} \text{and} P^{c(5)}$ to be log-normal distributed and $P^{c(6)}$ to be categorically distributed. Each evaluation run starts with an empty system.  We set the run lengths for each evaluation to be 300 days beyond a warm-up period of 300 days, which is discarded to avoid bias. Thirty initial design points are used to construct the Kriging model with two replications at each point, and an additional thirty points are selected sequentially based on the four approaches. The total evaluation budget is then $(30+ 30) \times 2 = 120$. 
As $f(\BFx, P^c)$ is unknown for this complex system, we instead use a long enough run lengths of $10^4$ days (as recommended in \cite{xie2014bayesian}) to estimate the system mean response, and then estimate the optimal solution $\BFx^* = \arg\min_{\BFx} f(\BFx, P^c)$ and optimal value $f(\BFx^*, P^c)$ based on the estimated mean response.  The estimated true optimal solution $\BFx^* $ is $[10,6,24]$ and the estimated optimal value $f(\BFx^*, P^c)$  is 1.248. 

\begin{figure}[!h]
	\centering
	\includegraphics[width=0.6\linewidth]{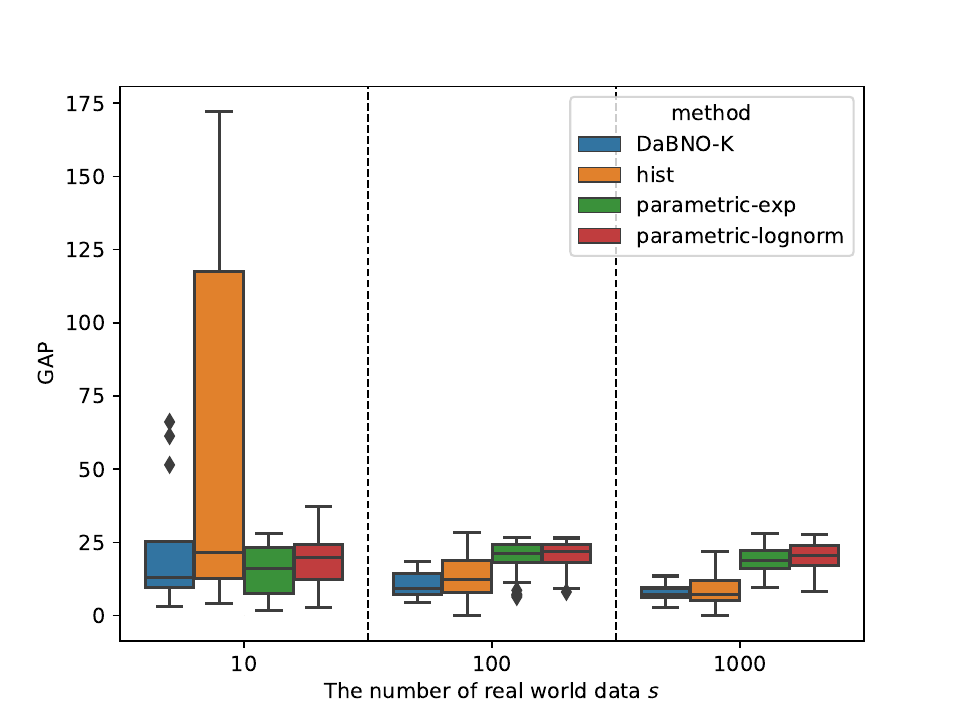}
	\caption{Inventory problem: boxplot of the GAP values for four approaches (DaBNO-K, hist, parametric-exponential, parametric-lognormal under different uncertainty levels of the distribution of $\BFu$.
	}
	\label{fig:inventory_gap}
\end{figure}

\begin{figure}
	\centering
	\includegraphics[width=0.6\linewidth]{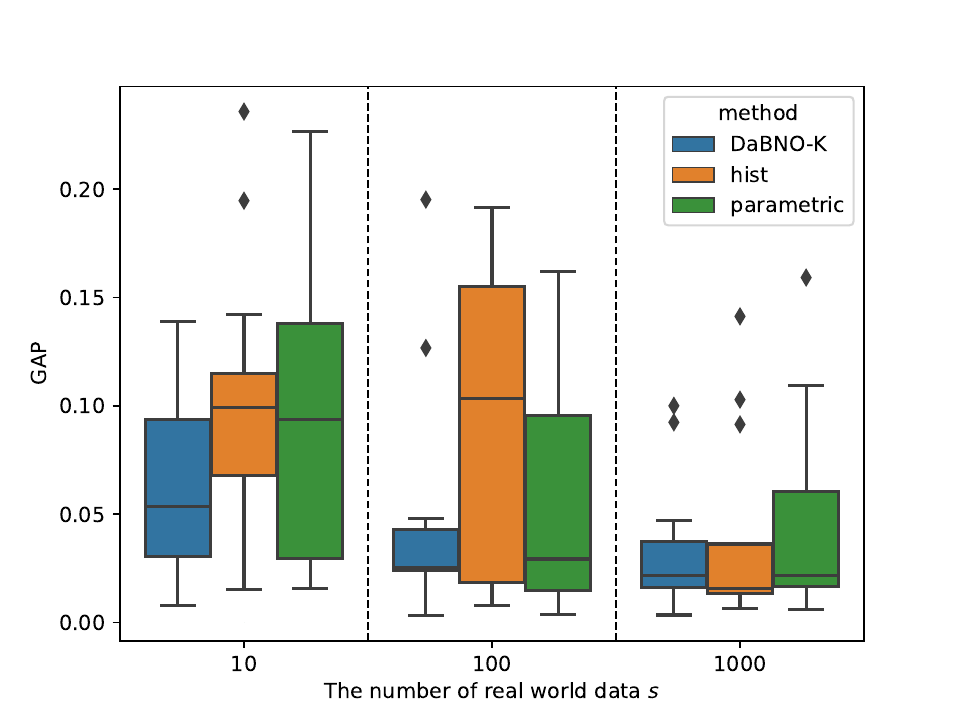}
	\caption{
		The critical care facility problem: the boxplot of the GAP values for three approaches (DaBNO-K, hist, parametric) under different uncertainty levels of the distribution of $\BFu$.}
	\label{fig:CCF_GAP}
\end{figure}

\paragraph{\textbf{Performance}}
For each algorithm, 100 trials are conducted to evaluate the performance and the GAP results are plotted in Figure \ref{fig:inventory_gap} and Figure \ref{fig:CCF_GAP}.

The results of the inventory problem are shown in Figure \ref{fig:inventory_gap}. Comparing DaBNO-K with the hist approach, it shows that the difference is not significant when $s=$ 1000 and the GAP values are close to zero with a small variance for both approaches. This result can be expected as the uncertainty of the distribution decreases when $s$ becomes large. When $s$ goes to $\infty$,  there is no uncertainty in the distribution of $\BFu$, and both the posterior $\BFpi(P|\BFD_s)$ and the empirical cumulative distributions $\hat{P}$ will be similar to the true distribution $P^c$. In this case, both approaches will reduce to the true problem (\ref{trueobjective}), and thus the minimizers from both approaches become close. 
In contrast, when $s=10$, the GAP values of the DaBNO-K approach are significantly better than those of the hist approach. 
When $s=100$,  the hist approach has a small GAP median but a relatively large variance. Comparing the performance of the DaBNO-K and the two parametric approaches, it shows that the DaBNO-K approach consistently outperforms the two parametric approaches across different uncertainty levels of uncertain distribution.  It also shows that the GAP values for the parametric-exponential approach and parametric-lognormal approach do not decrease as more real world data are available, which indicates that a wrong parametric distribution assumption will lead to very bad optimization results.

The results of the CCF problem are displayed in Figure \ref{fig:CCF_GAP}. It shows that when the number of real world data is small, i.e. $s = 10\text{ or }100$  the DaBNO-K approach outperforms hist and parametric approach with a lower mean and a relatively smaller variance. When $s = 1000$, DaBNO-K and hist approach perform quite similarly,  both outperform the parametric approach with a much smaller variance.

The good performance of our proposed algorithm comes at a cost of longer computational times, which is shown in Table \ref{tab:computation_time}. The main reason for the longer running time of DaBNO-K is that in the Kriging model, we need to calculate the Wasserstein distance between different distributions, which would require a lot of time both when we fit the Kriging model and when we make predictions. More efficient computation algorithms of the Wasserstein distance \citep{chizat2020faster} could be used to speed up the algorithm, and we leave this for future research. 
Although the computational time for our proposed approach is higher than those of the hist approach and parametric approaches, the differences are not too large and are reasonable for decisions required for longer periods. This also suggests that our method is well suited for expensive simulations (which are common in practice), where the evaluation time on the simulation model is much longer than the time required for selecting the subsequent evaluation point.
\begin{table}
	\caption{Computation time to run the whole algorithm. For the Griewantk, StybTang and inventory examples, the setting is 20 initial design points and  additional 40 evaluation points. The result for parametric is the average of parametric-exponential and parametric-lognormal.   For the critical care facility, the setting is 30 initial design points and additional 30 evaluation points. Unit: mins}
	\begin{center}
		\begin{tabular}{lccc}
			\hline
			\hline
			& DaBNO-K & hist & parametric \\
			\hline
			Griewank &82.5& 8.3&9.6\\
			StybTang &51.9& 7.4 &7.4\\
			Inventory & 85.4 & 13.8 & 13.6  \\
			Critical care facility       & 94.5 & 10.9 & 11.3 \\
			\hline
			\hline
		\end{tabular} 
	\end{center} 
	\label{tab:computation_time}
\end{table}

In summary, in comparison to the hist approach, the proposed DaBNO approach with the proposed DaBNO-K algorithm tends to be more robust to the uncertainty level of the uncertain distributions.
For the parametric approaches, it is clear that using the approaches that assume more information about the distribution family can sometimes bring about additional risks when the assumption involves errors. The additional risks can lead to sub-optimal or even wrong decisions. The proposed DaBNO-K that requires no information on the family of distributions is hence more robust.

\section{Conclusion} \label{sec3:conclusion}
In this paper, we proposed a data-drive  Bayesian nonparametric optimization (DaBNO) approach for the global optimization of black-box functions with uncertain parameters. We proposed to use the Dirichlet process to model the uncertainty in the distribution of the uncertain parameters and studied the asymptotic consistency and normality of the DaBNO to guarantee its theoretical performance. We also proposed a surrogate algorithm DaBNO-K based on EGO  to efficiently solve this DaBNO problem. 
We further studied the consistency properties of our DaBNO-K algorithm and tested its efficiency empirically with two synthetic functions, an inventory example and a critical care facility problem. We find that: 
1) The optimal solution/value of the DaBNO will converge the optimal solution/value of $f(\BFx,\BFP^c)$ as the amount of real-world data becomes large (uncertainty level decreases). 
2) The proposed algorithm DaBNO-K can converge to the optimal value/solution of $g(\BFx)$, which shows that our proposed algorithm can efficiently solve the DaBNO problem; and the true objective functions when uncertainty is low.
3) The proposed algorithm can obtain solutions closer to the optimal value of $f(\BFx,\BFP^c)$, compared with the approach ignoring the uncertainty in the distribution of $\BFu$. In addition, the DaBNO approach can avoid the risk of choosing a wrong parametric family and is more robust to the uncertainty level of the distribution of the uncertain parameters.

\section*{Acknowledgement}

Haowei Wang and Songhao Wang’s work is supported  by the National Natural Science Foundation (NNSF) of China under Grant 72101106 and the Shenzhen Science and Technology Program under Grant RCBS20210609103119020.
Szu Hui Ng’s work is supported in part by the Ministry of Education, Singapore (grant: R-266-000-149-114).

\bibliographystyle{informs2014} 
\bibliography{refhw} 



\ECSwitch


\ECHead{Proofs of Statements}
Throughout the proofs, we assume $\bm{u}_1,\cdots,\bm{u}_s$ are i.i.d random variable defined on the same probability space denoted by $(\Omega,\mathcal{F},P)$.
\section{Formulate the objective function as i.i.d summation}

\textbf{Proof of Lemma \ref{lemma:1}} 
Recall that $\bm{u}_1,...,\bm{u}_s\sim \bm{P}^c$. We model $\bm{P}^c$ with a prior $P\sim DP(\alpha,\bm{P}^0)$, where $DP(\alpha,\bm{P}^0)$ is a Dirichlet process. Denote the posterior distribution by $\bm{\pi}(\bm{P}|\bm{u}_1,...,\bm{u}_s)$. Then we have $\pi(P|\bm{u}_1,...,\bm{u}_s) \sim DP(\alpha+s, \frac{\alpha \bm{P}^0 + \sum_{j=1}^{s}\delta_{\bm{u}_j}}{\alpha+s})$.
Without loss of generality, we can assume $\overline{u}_i-\underline{u}_i=\delta_i$ for all $i=1,2,\cdots l $. For a fixed $m$,  considering the following partition  
\begin{equation*}
[\underline{u}_i,\overline{u}_i]=[\underline{u}_i,\underline{u}_i+\frac{\delta_i}{m}]\cup (\underline{u}_i+\frac{\delta_i}{m},\underline{u}_i+\frac{2\delta_i}{m}]\cup ...\cup (\overline{u}_i-\frac{\delta_i}{m},\overline{u}_i]=\cup_{k=1}^m B_k^{(i)}
\end{equation*}
\begin{equation*}
\begin{aligned}
\Omega&=\cup_{k_1,k_2,\cdots,k_l}  B_{k_1}^{(1)}\times B_{k_2}^{(2)}\times\cdots  \times B_{k_l}^{(l)}, \quad k_i=1,2,...,m, \text{for } i = 1, \cdots, l\\
&=\cup_{i=1}^{m^l} A_i (A_i = B_{k_1}^{(1)}\times B_{k_2}^{(2)}\times\cdots  \times B_{k_l}^{(l)}).
\end{aligned}
\end{equation*}
Now, denote the center of $A_i$ by by $\bm{c_i}$
and consider the piece-wise constant function 
\begin{equation*}
h_m(\bm{x},\bm{u})=\begin{cases}
h(\bm{x},\bm{c}_1),\quad \bm{u} \in A_1\\
h(\bm{x},\bm{c}_2),\quad \bm{u} \in A_2\\
...\\
h(\bm{x},\bm{c}_{m^{l}}),\quad \bm{u} \in A_{m^l}.
\end{cases}
\end{equation*}

For any distribution $\bm{P}$, we have
\begin{equation*}
	\mathbb{E}_{\bm{u} \sim \bm{P}} [h(\bm{x},\bm{u})]=\lim_{m\to +\infty} \sum_{i=1}^{m^l} h(\bm{x}, \bm{c}_i)\bm{P}(A_i)=\lim_{m\to +\infty} \mathbb{E}_{\bm{u} \sim \bm{P}} [h_m(\bm{x},\bm{u})]
\end{equation*}
Since $h(\bm{x},\bm{u})$ is Lipschitz, we know $\sup_{\bm{u}} |h_m(\bm{x},\bm{u})-h(\bm{x},\bm{u})|\leq \frac{L_1\delta}{m}$, where $\delta = \sqrt{\delta_1^2+ \delta_2^2 + \cdots +\delta_l^2}$, so
\begin{equation*}
\left|\sum_{i=1}^{m^l} h(\bm{x}, \bm{c}_i)\bm{P}(A_i)-\mathbb{E}_{\bm{P}} [h(\bm{x},\bm{u})]\right|=\left|\mathbb{E}_{\bm{P}} [h_m(\bm{x},\bm{u})]-\mathbb{E}_{\bm{P}} [h(\bm{x},\bm{u})]\right|\leq \frac{L_1\delta}{m}
\end{equation*}
As a result, if we put $Q_m(\bm{P})=\sum_{i=1}^{m^l} h(\bm{x},\bm{c}_i)\bm{P}(A_i)$, then 
$Q_m(\bm{P})\to \mathbb{E}_{\bm{P}} [h(\bm{x},\bm{u})]$ uniformly in $\bm{P}$. We now compute
\begin{equation}
\mathbb{E}_{\bm{\pi}}[\mathbb{E}_{\bm{P}} [h(\bm{x},\bm{u})]]=\mathbb{E}_{\bm{\pi}} \left[\lim_m Q_m(\bm{P})\right]=\lim_m \mathbb{E}_{\bm{\pi}} [Q_m(\bm{P})]. \label{limexp}
\end{equation}
We can exchange the order of limit and expectation in (\ref{limexp}) because $Q_m(\bm{P})\to \mathbb{E}_{\bm{P}} [h(\bm{x},\bm{u})]$ uniformly in $\bm{P}$.
For any fixed $m$, we can compute
\begin{equation*}
\begin{aligned}
\mathbb{E}_{\pi} [Q_m(\bm{P})]&=\sum_{i=1}^m h(\bm{x},\bm{c}_i)\mathbb{E}_{\pi}(\bm{P}(A_i))=\sum_{i=1}^m h(\bm{x},\bm{c}_i)\left( \frac{\alpha \bm{P}^0(A_i)}{\alpha+s}+\frac{1}{\alpha+s}\sum_{j=1}^s \delta_{\bm{u}_j}(A_i) \right)
\end{aligned}
\end{equation*}
The first term of the above equation converges to $\mathbb{E}_{\bm{P}^0}\frac{\alpha}{\alpha+s}h(\bm{x},\bm{u})$ by definition.
We now consider the second term. 
For fixed $s$, when $m$ tends to  infinity,
\begin{equation*}
\resizebox{\hsize}{!}{$
	\begin{aligned}
	\lim_{m\to+\infty}\sum_{i=1}^m\left(h(\bm{x},\bm{c}_i)\times \frac{1}{\alpha+s}\sum_{j=1}^s \delta_{\bm{u}_j}(A_i)\right)
	&=\frac{1}{\alpha+s}\sum_{j=1}^s \lim_{m\to+\infty}\sum_{i=1}^m h(\bm{x},\bm{c}_i)\delta_{\bm{u}_j}(A_i)\\
	&=\frac{1}{\alpha+s}\sum_{j=1}^s \lim_{m\to +\infty} h(\bm{x},\bm{c}_{i'}) ( \mbox{such that } \bm{u}_j \in A_{i'}, {i'} \in {1, \dots, m^l})\\
	&=\frac{1}{\alpha+s}\sum_{j=1}^s h(\bm{x},\bm{u}_j).
	\end{aligned}
	$}
\end{equation*}
The second equality holds because the partitions $A_i, i = 1, \cdots, m^l $ are disjoint,  implying that $\bm{u}_j$ can only belong to one of the partitions: now suppose $\bm{u}_j \in A_{i'}, {i'} \in [1, \dots, m^l])$.
Combining all this, we have $\mathbb{E}_{\bm{\pi}}[\mathbb{E}_{\bm{P}} [h(\bm{x},\bm{u})]]=\mathbb{E}_{\bm{P}^0}\frac{\alpha}{\alpha+s}h(\bm{x},\bm{u})+\frac{1}{\alpha+s}\sum_{j=1}^s h(\bm{x},\bm{u}_j).$

\section{Proof of consistency (Theorem \ref{theorem:1}).}
\textbf{Proof of $(i)$:}	According to Lemma \ref{lemma:1}, by adding and subtracting $\frac{1}{s}\sum_{j=1}^sh(\bm{x},\bm{u}_j)$, we have
\begin{equation*}
\begin{aligned}
&\mathbb{E}_{\bm{\pi}}[\mathbb{E}_{\bm{P}}[ h(\bm{x},\bm{u})]]-\mathbb{E}_{\bm{P}^c}[h(\bm{x},\bm{u})]\\
&= \mathbb{E}_{\bm{P}^0}\frac{\alpha}{\alpha+s}h(\bm{x},\bm{u})+\frac{1}{\alpha+s}\sum_{j=1}^s h(\bm{x},\bm{u}_j)-\mathbb{E}_{\bm{P}^c}[h(\bm{x},\bm{u})]\\
&=\mathbb{E}_{\bm{P}^0}\frac{\alpha}{\alpha+s}h(\bm{x},\bm{u})+\frac{1}{\alpha+s}\sum_{j=1}^s h(\bm{x},\bm{u}_j)-\frac{1}{s}\sum_{j=1}^sh(\bm{x},\bm{u}_j)+\frac{1}{s}\sum_{j=1}^s\left(h(\bm{x},\bm{u}_j)- \mathbb{E}_{\bm{P}^c}[h(\bm{x},\bm{u})]\right)\\
&=\frac{\alpha}{\alpha+s}\mathbb{E}_{\bm{P}^0} h(\bm{x},\bm{u})-\frac{\alpha}{(s+\alpha){s}}\sum_{j=1}^s h(\bm{x},\bm{u}_j)+\frac{1}{s}\sum_{j=1}^s \left(h(\bm{x},\bm{u}_j)- \mathbb{E}_{\bm{P}^c}[h(\bm{x},\bm{u})]\right) 
\end{aligned}
\end{equation*}
The first two terms of the above equation are of order $o(1)$.
We then apply the strong law of large  numbers to the last term of the above equation to conclude.

\textbf{Proof of  $(ii)$:} Since for any $\bm{u}$, $|h(\bm{x}_1,\bm{u})-h(\bm{x}_2,\bm{u})|\leq L_2 ||\bm{x}_1-\bm{x}_2||$ according to Assumption \ref{assum}, the class of functions $\mathcal{F}=\{h(\bm{x},\bm{u})| \bm{x}\in \mathbb{R} \}$ is Donsker, from Example 19.6 of \cite{vandervarrt}, and thus $|\mathbb{E}_{\bm{P}^c}[h(\bm{x},\bm{u})]-\frac{1}{s}\sum_{j=1}^s h(\bm{x},\bm{u}_j)|\to_{a.s.} 0$ uniformly by the property of Donsker class functions \citep{vandervarrt}. 
Therefore, 
\begin{equation}
\sup_{\bm{x}}|\mathbb{E}_{\bm{P}^c}[h(\bm{x},\bm{u})]-\frac{1}{s}\sum_{j=1}^s h(\bm{x},\bm{u}_j)| \to_{a.s.} 0.
\label{equ:con2}
\end{equation}
Next, we know
\begin{equation*}
\begin{aligned}
\sup_{\bm{x}} |\mathbb{E}_{\bm{\pi}}[\mathbb{E}_{\bm{P}} [h(\bm{x},\bm{u})]] -\frac{1}{s}\sum_{j=1}^s h(\bm{x},\bm{u}_j)|&= \sup_{\bm{x}}|\frac{\alpha}{\alpha+s}\mathbb{E}_{\bm{P}^0} h(\bm{x},\bm{u})-\frac{\alpha}{(s+\alpha)s}\sum_{j=1}^s h(\bm{x},\bm{u}_j)| \\
& \leq \frac{\alpha}{\alpha+s} \sup_{\bm{x}} \left( |\mathbb{E}_{\bm{P}^0} h(\bm{x},\bm{u})| + |\frac{1}{s}\sum_{j=1}^s h(\bm{x},\bm{u}_j)|  \right)\\
&\leq \frac{2\alpha}{\alpha+s} \sup_{\bm{x},\bm{u}} |h(\bm{x},\bm{u})|=o(1).
\end{aligned}
\end{equation*}
Combine the above display and \eqref{equ:con2} and use triangle inequality to establish that
\begin{align}\label{eq:uniform as}
\sup_{\bm{x}}|\mathbb{E}_{\bm{P}^c}[h(\bm{x},\bm{u})]-\mathbb{E}_{\bm{\pi}}[\mathbb{E}_{\bm{P}} [h(\bm{x},\bm{u})]]|& \leq \sup_{\bm{x}} |\mathbb{E}_{\bm{\pi}}[\mathbb{E}_{\bm{P}} [h(\bm{x},\bm{u})]] -\frac{1}{s}\sum_{j=1}^s h(\bm{x},\bm{u}_j)|\notag\\
&\quad+\sup_{\bm{x}}|\mathbb{E}_{\bm{P}^c}[h(\bm{x},\bm{u})]-\frac{1}{s}\sum_{j=1}^s h(\bm{x},\bm{u}_j)|  \to_{a.s.} 0.
\end{align}
Denote $g(\bm{x})\triangleq   \mathbb{E}_{\bm{\pi}}[\mathbb{E}_{\bm{P}} [h(\bm{x},\bm{u})]] = \mathbb{E}_{\bm{P}^0}\frac{\alpha}{\alpha+s}h(\bm{x},\bm{u})+  \frac{1}{\alpha+s}\sum_{i=1}^s h(\bm{x},\bm{u}_i)$ and $f(x)\triangleq \mathbb{E}_{\bm{P}^c}[h(\bm{x},\bm{u})]$.
Also denote the minimizer of $g$ as $\bm{x}_s$, the set of minimizer of $f$ as $S$. Here, we actually prove a stronger result where we no-longer require $S$ contains only one element but assume that $S$ is a bounded compact set (possibly contains multiple elements). 
It can be shown that $\mbox{dist}(\bm{x}_s,S)\to_{a.s.} 0$, where 
\begin{align*}
\mbox{dist}(\bm{x},S)\triangleq min_{\bm{y}\in S}\| \bm{x}-\bm{y} \|.
\end{align*}
Every $\omega\in \tilde{\Omega}$ corresponds to a sequence of realization of $\{\bm{u}_i \}_{i=1}^{\infty}$.
To prove the almost sure convergence, we suffice to show for almost every $\omega\in\tilde{\Omega}$,  $\mbox{dist}(\bm{x}_s,S)\to 0 $ (a deterministic result).
By \eqref{eq:uniform as} we have for almost every  $\omega\in\tilde{\Omega}$ the following holds
\begin{equation}\label{eq: as conv}
\sup_{\bm{x}} \left|\mathbb{E}_{\bm{P}^0}\frac{\alpha}{\alpha+s}h(\bm{x},\bm{u})+  \frac{1}{\alpha+s}\sum_{i=1}^s h(\bm{x},\bm{u}_i)-\mathbb{E}_{\bm{P}^c}[h(\bm{x},\bm{u})] \right|\to 0.
\end{equation}
We claim that for every $\omega$ such that the above display holds, $\mbox{dist}(\bm{x}_s,S)\to 0 $.
Suppose otherwise this is not true, and we have a sequence $\{\bm{x}_k \}$ such that $\mbox{dist}(\bm{x}_k,S)\nrightarrow 0$.
Since all the $\bm{x}_k$ are contained in a compact set $\mathscr{X}$,  by passing to a subsequence if necessary, we can assume that $\mbox{dist}(\bm{x}_k,S)\geq \epsilon>0$, and that $\bm{x}_k\to \bm{x}_0\in \mathscr{X}$. It follows that $\bm{x}_0\notin S$ and $f(\bm{x}_0)>\min_{\bm{x}} f(\bm{x}).$
Furthermore, we can write
\begin{equation*}
g(\bm{x}_s)-f(\bm{x}_0)=g(\bm{x}_s)-f(\bm{x}_s)+f(\bm{x}_s)-f(\bm{x}_0).
\end{equation*}
The first term on the RHS converges to zero  from \eqref{eq: as conv}.
The second term converges to zero from our assumption on $\bm{x}_0$.
We conclude that for every $\omega\in \tilde{\Omega}$ such that $ \sup_{\bm{x}} |g(\bm{x})-f(\bm{x})|\to 0$, we must have $g(\bm{x}_s)-f(\bm{x}_0)\to 0$.
However, we also have $|g(\bm{x}_s)-\min_{\bm{x}} f(\bm{x})|\leq  \sup_{\bm{x}} |g(\bm{x})-f(\bm{x})|\to 0$ for such $\omega$.
This contradicts with the fact $f(\bm{x}_0)>\min_x f(\bm{x}).$
This means that our initial assumption is wrong and so $\mbox{dist}(\bm{x}_s,S)\to 0 $ for every $\omega$ such that $ \sup_{\bm{x}} |g(\bm{x})-f(\bm{x})|\to 0$.
Such $\omega$ is almost everywhere.

\textbf{Proof of $(iii)$:}
The almost sure convergence of optimal value follows directly from the above part and \eqref{eq:uniform as}.
Specifically, recall the notation in the previous part $g(\bm{x})=   \mathbb{E}_{\bm{\pi}}[\mathbb{E}_{\bm{P}} [h(\bm{x},\bm{u})]] = \mathbb{E}_{\bm{P}^0}\frac{\alpha}{\alpha+s}h(\bm{x},\bm{u})+  \frac{1}{\alpha+s}\sum_{i=1}^s h(\bm{x},\bm{u}_i)$ and $f(x)=\mathbb{E}_{\bm{P}^c}[h(\bm{x},\bm{u})]$.
We write
\begin{align*}
g(\bm{x}_s)-f(\bm{x}^*)=g(\bm{x}_s)-f(\bm{x}_s)+f(\bm{x}_s)-f(\bm{x}^*).
\end{align*}
The first term on the RHS converges to zero almost surely from \eqref{eq:uniform as} and the second term is from the continuity of $f$ and the almost sure convergence of $x_s$.

\section{Proof of asymptotic normality (Theorem \ref{theorem:2})}
\textbf{Proof of $(i)$}:	According to Lemma \ref{lemma:1},  by adding and subtracting $\frac{1}{s}\sum_{j=1}^sh(\bm{x},\bm{u}_j)$, we have
\begin{equation*}
\resizebox{1\hsize}{!}{$
	\begin{aligned}
	&\sqrt{s}\left(\mathbb{E}_{\bm{\pi}}[\mathbb{E}_{\bm{P}}[ h(\bm{x},\bm{u})]]-\mathbb{E}_{\bm{P}^c}[h(\bm{x},\bm{u})]\right)\\
	&=\sqrt{s}\left( \mathbb{E}_{\bm{P}^0}\frac{\alpha}{\alpha+s}h(\bm{x},\bm{u})+\frac{1}{\alpha+s}\sum_{j=1}^s h(\bm{x},\bm{u}_j)-\mathbb{E}_{\bm{P}^c}[h(\bm{x},\bm{u})] \right)\\
	&=\sqrt{s}\left( \mathbb{E}_{\bm{P}^0}\frac{\alpha}{\alpha+s}h(\bm{x},\bm{u})+\frac{1}{\alpha+s}\sum_{j=1}^s h(\bm{x},\bm{u}_j)-\frac{1}{s}\sum_{j=1}^sh(\bm{x},\bm{u}_j)+\frac{1}{s}\sum_{j=1}^s\left(h(\bm{x},\bm{u}_j)- \mathbb{E}_{\bm{P}^c}[h(\bm{x},\bm{u})]\right)  \right)\\
	&=\frac{\sqrt{s}\alpha}{\alpha+s}\mathbb{E}_{\bm{P}^0} h(\bm{x},\bm{u})-\frac{\alpha}{(s+\alpha)\sqrt{s}}\sum_{j=1}^s h(\bm{x},\bm{u}_j)+\frac{1}{\sqrt{s}}\sum_{j=1}^s \left(h(\bm{x},\bm{u}_j)- \mathbb{E}_{\bm{P}^c}[h(\bm{x},\bm{u})]\right) 
	\end{aligned}
	$}
\end{equation*}
The first two terms of the above equation are of order $o_p(1)$, which do not affect the asymptotic distribution of $\sqrt{s}\left(\mathbb{E}_{\bm{\pi}}[\mathbb{E}_{\bm{P}}[ h(\bm{x},\bm{u})]]-\mathbb{E}_{\bm{P}^c}[h(\bm{x},\bm{u})]\right) $.
We then apply the central limit theorem to the last term of the above equation to conclude.

\textbf{Proof of $(ii)$: }
Since $\bm{x}_s$ minimize $\mathbb{E}_{\bm{\pi}}[\mathbb{E}_{\bm{P}}[ h(\bm{x},\bm{u})]]$, 
$$\mathbb{E}_{\bm{P}^0}\frac{\alpha}{\alpha+s}\partial_{\bm{x}} h(\bm{x}_s,\bm{u})+\frac{1}{\alpha+s}\sum_{j=1}^s \partial_{\bm{x}} h(\bm{x}_s,\bm{u}_j)=0. $$
Using the above equation and the fact $\frac{1}{s}\sum_{j=1}^s \partial_{\bm{x}} h(\bm{x}_s,\bm{u}_j)=\frac{1}{\alpha+s}\sum_{j=1}^s \partial_{\bm{x}} h(\bm{x}_s,\bm{u}_j)+\frac{\alpha}{s(s+\alpha)}\sum_{j=1}^s \partial_{\bm{x}} h(\bm{x}_s,\bm{u}_j) $, the following equality is immediate:
\begin{equation*}
\begin{aligned}
\frac{1}{s}\sum_{j=1}^s \partial_{\bm{x}} h(\bm{x}_s,\bm{u}_j)&= \mathbb{E}_{\bm{P}^0}\frac{\alpha}{s+\alpha}\partial_{\bm{x}} h(\bm{x}_s,\bm{u})+\frac{1}{\alpha+s}\sum_{j=1}^s \partial_{\bm{x}} h(\bm{x}_s,\bm{u}_j)-\\
&\left( \mathbb{E}_{\bm{P}^0}\frac{\alpha}{s+\alpha}\partial_{\bm{x}} h(\bm{x}_s,\bm{u})-\frac{\alpha}{s(s+\alpha)}\sum_{j=1}^s \partial_{\bm{x}} h(\bm{x}_s,\bm{u}_j)  \right)=o_p(\frac{1}{\sqrt{s}})
\end{aligned}
\end{equation*}
The last equality holds due to the fact that $ \mathbb{E}_{\bm{P}^0}\frac{\alpha}{s+\alpha}\partial_{\bm{x}} h(\bm{x}_s,\bm{u})-\frac{\alpha}{s(s+\alpha)}\sum_{j=1}^s \partial_{\bm{x}} h(\bm{x}_s,\bm{u}_j) =O_p(\frac{1}{s}) = o_p(\frac{1}{\sqrt{s}})$.
Since $\partial_{\bm{x}} h(\bm{x},\bm{u})$ is Lipschitz in $\bm{x}$ (Assumption \ref{assum:1}$(i)$), $\bm{x}_s \to_p \bm{x}^*$,  and $\frac{1}{s}\sum_{j=1}^s \partial_{\bm{x}} h(\bm{x}_s,\bm{u}_j) = o_p(\frac{1}{\sqrt{s}})$, we have verified the conditions  in Theorem 5.21 of \cite{vandervarrt}, and so $\sqrt{s}(\bm{x}_s-\bm{x}^*)$ is asymptotically normal distributed with mean zero and asymptotic covariance matrix $$\Sigma=V_{{\bm{x}}^*}^{-1}\mathbb{E}\frac{\partial}{\partial {\bm{x}}}h({\bm{x}^*},\bm{u})\frac{\partial}{\partial {\bm{x}}}h({\bm{x}^*},\bm{u})^T (V_{{\bm{x}}^*}^{-1})^{T}.$$  from Theorem 5.21 of \cite{vandervarrt}.

\textbf{Proof of $(iii)$:}
Since for any $\bm{u}$, $|h(\bm{x}_1,\bm{u})-h(\bm{x}_2,\bm{u})|\leq L_2 ||\bm{x}_1-\bm{x}_2||$ from Assumption \ref{assum}, the class of functions $\mathcal{F}=\{h(\bm{x},\bm{u})| x\in \mathbb{R} \}$ is Donsker, from Example 19.6 of \cite{vandervarrt}. 
Also since $ \bm{x}_s \to_p \bm{x}^*$, 
we can use Lemma 19.24 of \cite{vandervarrt} to see
\begin{equation*}
	\begin{aligned}
	\frac{1}{\sqrt{s}}\sum_{j=1}^s h(\bm{x}_s,\bm{u}_j)-\frac{1}{\sqrt{s}}\sum_{j=1}^s h(\bm{x}^*,\bm{u}_j)=\sqrt{s}\left( \mathbb{E}_{\bm{P}^c}[h(\bm{x}_s,\bm{u})]-\mathbb{E}_{\bm{P}^c}[h(\bm{x}^*,\bm{u})] \right)+o_p(1).
	\end{aligned}
\end{equation*}
Using the Taylor expansion to deal with  the right-hand side of the above equation, and using the fact $(\mathbb{E}_{\bm{P}^c} [h(\bm{x},\bm{u})])'_{\bm{x} =\bm{x}^*}=0$, we obtain
\begin{equation}
\frac{1}{\sqrt{s}}\sum_{j=1}^s h(\bm{x}_s,\bm{u}_j)-\frac{1}{\sqrt{s}}\sum_{j=1}^s h(\bm{x}^*,\bm{u}_j)=o_p(\sqrt{s}|\bm{x}_s-\bm{x}^*|). 
\label{equ:consist}
\end{equation}
Thus
\begin{equation}\label{eq: consist h(xn h(x*))}
\frac{1}{s}\sum_{j=1}^sh(\bm{x}_s,\bm{u}_j)-\frac{1}{s}\sum_{j=1}^sh(\bm{x}^*,\bm{u}_j)=o_p(1).
\end{equation}	
By using Lemma \ref{lemma:1} and telescoping with $\frac{1}{s}\sum_{j=1}^sh(\bm{x}_s,\bm{u}_j)$, we have
\begin{equation*}
\resizebox{1\hsize}{!}{$
	\begin{aligned}
	&\sqrt{s}\left(\mathbb{E}_{\bm{\pi}}[\mathbb{E}_{\bm{P}}[ h(\bm{x}_s,\bm{u})]]-\mathbb{E}_{\bm{P}^c}[h(\bm{x}^*,\bm{u})]\right)\\
	&=\sqrt{s}\left(\mathbb{E}_{\bm{P}^0}\frac{\alpha}{\alpha+s}h(\bm{x}_s,\bm{u})+\frac{1}{\alpha+s}\sum_{j=1}^s h(\bm{x}_s,\bm{u}_j)-\mathbb{E}_{\bm{P}^c}[h(\bm{x}^*,\bm{u})]\right)\\
	&=\sqrt{s}\left(\mathbb{E}_{\bm{P}^0}\frac{\alpha}{\alpha+s}h(\bm{x}_s,\bm{u})+\frac{1}{\alpha+s}\sum_{j=1}^s h(\bm{x}_s,\bm{u}_j)-\frac{1}{s}\sum_{j=1}^sh(\bm{x}_s,\bm{u}_j)+\frac{1}{s}\sum_{j=1}^sh(\bm{x}_s,\bm{u}_j)- \mathbb{E}_{\bm{P}^c}[h(\bm{x}^*,\bm{u})]\right)\\
	&= \mathbb{E}_{\bm{P}^0}\frac{\sqrt{s}\alpha}{\alpha+s}h(\bm{x}_s,\bm{u}) - \frac{\alpha}{(\alpha+s)\sqrt{s}}\sum_{j=1}^s h(\bm{x}_s,\bm{u}_j) + \frac{1}{\sqrt{s}}\sum_{j=1}^sh(\bm{x}_s,\bm{u}_j)- \sqrt{s} \mathbb{E}_{\bm{P}^c} [h(\bm{x}^*,\bm{u})]\\
	&  = \mathbb{E}_{\bm{P}^0}\frac{\sqrt{s}\alpha}{\alpha+s}h(\bm{x}_s,\bm{u}) - \frac{\alpha}{(\alpha+s)\sqrt{s}}\sum_{j=1}^s h(\bm{x}_s,\bm{u}_j)   + \frac{1}{\sqrt{s}}\sum_{j=1}^sh(\bm{x}_s,\bm{u}_j) -\frac{1}{\sqrt{s}}\sum_{j=1}^sh(\bm{x}^*,\bm{u}_j)  \\
	&\quad +\frac{1}{\sqrt{s}}\sum_{j=1}^sh(\bm{x}^*,\bm{u}_j) - \sqrt{s}\mathbb{E}_{\bm{P}^c} [h(\bm{x}^*,\bm{u})] 
	\end{aligned}
	$}
\end{equation*}
The sum of the first two terms of the last equality is of order $o_p(1)$ since $h$ is bounded. We can then use the fact that $\frac{1}{\sqrt{s}}\sum_{j=1}^sh(\bm{x}_s,\bm{u}_j)-\frac{1}{\sqrt{s}}\sum_{j=1}^sh(\bm{x}^*,\bm{u}_j)=o_p(1)$ to reformulate the last equality of the above equation as
\begin{equation*}
\begin{aligned}
&\sqrt{s}\left(\mathbb{E}_{\bm{\pi}}[\mathbb{E}_{\bm{P}}[ h(\bm{x}_s,\bm{u})]]-\mathbb{E}_{\bm{P}^c}[h(\bm{x}^*,\bm{u})]\right) =o_p(1) + \frac{1}{\sqrt{s}}\sum_{j=1}^s\left(h(\bm{x}^*,\bm{u}_j)-\mathbb{E}_{\bm{P}^c} [h(\bm{x}^*,\bm{u})]\right)
\end{aligned}
\end{equation*}
According to the central limit theorem, $\frac{1}{\sqrt{s}}\sum_{j=1}^s\left(h(\bm{x}^*,\bm{u}_j)-\mathbb{E}_{\bm{P}^c} [h(\bm{x}^*,\bm{u})]\right)$ is asymptotically normally distributed with mean zero and variance $\text{Var}(h(\bm{x}^*,\bm{u}))$.

\section{Proof of the convergence of algorithm}
\textbf{Proof of Theorem \ref{thm:main theorem}}
Recall that the proof decomposes into two parts.
First look at the first step, since 
\begin{equation*}
\left|\mu_n(\bm{x})-g(\bm{x})\right|\leq \left|\mu_n(\bm{x})-\mathbb{E}_{\bm{P} \sim \bm{\pi}}[F_n(\bm{x},\bm{P})]\right|+\left|\mathbb{E}_{\bm{P} \sim \bm{\pi}}[F_n(\bm{x},\bm{P})]-g(\bm{x})   \right|,
\end{equation*}
from the triangle inequality,
\begin{equation*} 
\begin{aligned}
&\mathbb{P}\left\lbrace\left|\mu_n(\bm{x})-g(\bm{x})\right| \geq 2\epsilon \right\rbrace \\&\leq \mathbb{P}\left\lbrace\left|\mu_n(\bm{x})-\mathbb{E}_{\bm{P} \sim \bm{\pi}}F_n(\bm{x},\bm{P})\right|+\left|\mathbb{E}_{\bm{P} \sim \bm{\pi}}F_n(\bm{x},\bm{P})-g(\bm{x})   \right|\geq 2\epsilon \right\rbrace\\
&\leq \mathbb{P}\left\lbrace \left|\mu_n(\bm{x})-\mathbb{E}_{\bm{P} \sim \bm{\pi}}F_n(\bm{x},\bm{P})\right|\geq \epsilon \right\rbrace+\mathbb{P}\left\lbrace \left|\mathbb{E}_{\bm{P} \sim \bm{\pi}}F_n(\bm{x},\bm{P})-g(\bm{x}) \right|\geq \epsilon \right\rbrace.
\end{aligned}
\end{equation*}

Thus, to establish Equation \eqref{eq: mu_n to g}, it suffices to establish
\begin{equation*}
\mu_n(\bm{x})-\mathbb{E}_{\bm{P} \sim \bm{\pi}}F_n(\bm{x},\bm{P})\to_p 0,\quad \mbox{and} \quad \mathbb{E}_{\bm{P} \sim \bm{\pi}}F_n(\bm{x},\bm{P})-g(\bm{x}) \to_p 0
\end{equation*}
Now, suppose we could directly observe $f(\bm{x},\bm{P})$ and construct a deterministic Kriging model with conditional mean $\tilde{m}_n(\bm{x},\bm{P})$ and conditional variance $\tilde{k}_n\left( (\bm{x},\bm{P}),(\bm{x},\bm{P}) \right)$. 
We have
\begin{equation*}
\resizebox{1\hsize}{!}{$
	\begin{aligned}
	\left|F_n(\bm{x},\bm{P})-f(\bm{x},\bm{P})\right|&\leq \underbrace{\left|F_n(\bm{x},\bm{P})-m_n(\bm{x},\bm{P})  \right|}_{A}+\underbrace{\left| m_n(\bm{x},\bm{P})-\tilde{m}_n(\bm{x},\bm{P})\right|}_{B}+\underbrace{\left| \tilde{m}_n(\bm{x},\bm{P})-f(\bm{x},\bm{P}) \right|}_{C}
	\end{aligned}
	$}
\end{equation*}
Direct calculations show that 
\begin{equation*}
\tilde{m}_n(\bm{x},\bm{P})-m_n(\bm{x},\bm{P})\to_p 0 \quad \mbox{and} \quad \tilde{k}_n\left( (\bm{x},\bm{P}),(\bm{x},\bm{P}) \right)-k_n\left( (\bm{x},\bm{P}),(\bm{x},\bm{P}) \right)\to_p 0
\end{equation*}
uniformly in $\bm{x}$ as the number of replications $r$ tends to infinity (Lemma 5 of \cite{pedrielli2020extended}). 
And so part $B$ tends to zero uniformly.
According to Propositions 3.3, 3.4 and 3.5 of \cite{AMS18}, we have
$
\tilde{m}_n(\bm{x},\bm{P})-f(\bm{x},\bm{P})\to_p 0, \quad \mbox{and} \quad k_n\left( (\bm{x},\bm{P}),(\bm{x},\bm{P}) \right)\to_p 0,
$
uniformly for all $(\bm{x},\bm{P})$.
And so part $C$ tends to zero as $s$ tends to infinity.
For part $A$, the expectation of $F_n(\bm{x},\bm{P})$ is $m_n(\bm{x},\bm{P})$ and the marginal variance $k_n$ tends to zero uniformly. Thus part $A$ also converges to zero in probability uniformly for $(\bm{x},\bm{P})$ from Chebyshev's Inequality:
\begin{equation*}
\mathbb{P}(\left|F_n(\bm{x},\bm{P})-m_n(\bm{x},\bm{P})  \right|\geq \epsilon)\leq \frac{k_n((\bm{x},\bm{P}),(\bm{x},\bm{P}))}{\epsilon^2}.
\end{equation*}
Combining all that we know,
\begin{equation*}
\left|F_n(\bm{x},\bm{P})-f(\bm{x},\bm{P})\right|\to_p 0 
\end{equation*}
uniformly for all $(\bm{x},\bm{P})$.
And so
\begin{equation}
\begin{aligned}
\sup_{\bm{x}}\left|\mathbb{E}_{\bm{P} \sim \bm{\pi}}F_n(\bm{x},\bm{P})-g(\bm{x})\right|
&=\sup_{\bm{x}}\left|\mathbb{E}_{\bm{P} \sim \bm{\pi}}\left(F_n(\bm{x},\bm{P})-f(\bm{x},\bm{P})\right)\right|\\
&\leq \sup_{\bm{x},\bm{P}}\left|\left(F_n(\bm{x},\bm{P})-f(\bm{x},\bm{P})\right)\right|\to_p 0,
\end{aligned}
\label{eq:errorpart2}
\end{equation}
i.e., $ \mathbb{E}_{\bm{P} \sim \bm{\pi}}F_n(\bm{x},\bm{P})-g(\bm{x}) \to_p 0$ uniformly in $\bm{x}$. In addition,
\begin{equation*}
\begin{aligned}
\left| \mu_n(\bm{x})-\mathbb{E}_{\bm{P} \sim \bm{\pi}}F_n(\bm{x},\bm{P}) \right| &=\left| \mu_n(\bm{x})-\mathbb{E}_{\bm{P} \sim \bm{\pi}}m_n(\bm{x},\bm{P})+\mathbb{E}_{\bm{P} \sim \bm{\pi}}m_n(\bm{x},\bm{P})-\mathbb{E}_{\bm{P} \sim \bm{\pi}}F_n(\bm{x},\bm{P}) \right|\\
&\leq \left| \mu_n(\bm{x})-\mathbb{E}_{\bm{P} \sim \bm{\pi}}m_n(\bm{x},\bm{P})\right|+\left|\mathbb{E}_{\bm{P} \sim \bm{\pi}}m_n(\bm{x},\bm{P})-\mathbb{E}_{\bm{P} \sim \bm{\pi}}F_n(\bm{x},\bm{P}) \right|
\end{aligned}
\end{equation*}
Recall
\begin{equation*}
\mu_n(\bm{x})=\frac{1}{N_{mc}}\sum_{i=1}^{N_{mc}} m_n(\bm{x}, \bm{P}_i).
\end{equation*}
By definition, $m_n(\bm{x},\bm{P})$ is uniformly bounded for all $\bm{x}$ and $\bm{P}$ and so Chebyshev's Inequality implies
\begin{equation*}
\mathbb{P}\left\lbrace \left| \mu_n(\bm{x})-\mathbb{E}_{\bm{P} \sim \bm{\pi}}m_n(\bm{x},\bm{P})\right|\geq \epsilon \right\rbrace\leq \frac{\mathbb{E}_{\bm{P} \sim \bm{\pi}} m_n^2(\bm{x},\bm{P}) - (\mathbb{E}_{\bm{P} \sim \bm{\pi}} m_n(\bm{x},\bm{P}))^2}{\epsilon^2N_{mc}}.
\end{equation*}
And so $ \mu_n(\bm{x})-\mathbb{E}_{\bm{P} \sim \bm{\pi}}m_n(\bm{x},\bm{P})\to_p 0$ uniformly in $\bm{x}$ as $N_{mc} \to \infty$.
On the other hand, 
\begin{equation*}
\sup_{\bm{x}} \left|\mathbb{E}_{\bm{P} \sim \bm{\pi}}m_n(\bm{x},\bm{P})-\mathbb{E}_{\bm{P} \sim \bm{\pi}}F_n(\bm{x},\bm{P})\right| \leq \sup_{\bm{x},\bm{P}} \left|  m_n(\bm{x},\bm{P})-F_n(\bm{x},\bm{P}) \right| 
\end{equation*}
The right-hand side converges to zero in probability uniformly and so $\mathbb{E}_{\bm{P} \sim \bm{\pi}}m_n(\bm{x},\bm{P})-\mathbb{E}_{\bm{P} \sim \bm{\pi}}F_n(\bm{x},\bm{P})\to_p 0 $.
So $\mu_n(\bm{x})-\mathbb{E}_{\bm{P} \sim \bm{\pi}}F_n(\bm{x},\bm{P})\to_p 0$ uniformly in $\bm{x}$.
In view of Equation (\ref{eq:errorpart2}), we have completed the proof of $\left|\mu_n(\bm{x})-g(\bm{x})\right|\to_p 0.$

We then move on to the second step of the proof, i.e., to establish the density of the points we visit.
The proof follows similar procedures as Lemma 1 of \cite{locatelli1997bayesian}, albeit in a multiple dimensional setting as in \cite{wangsh2019multi}.
We consider the following stopping rule in the algorithm: we continue the algorithm  while the maximum expected improvement $\mathbb{EI}_T(\bm{x}_{n+1},\bm{P}_{n+1})$ is above some pre-determined constant $c$ (which is allowed to be arbitrary small in this setting). 
Recall the definition of $ \text{EI}_T(\bm{x}_{n+1},\bm{P}_{n+1}):$
\begin{equation*}
\text{EI}_T(\bm{x}_{n+1},\bm{P}_{n+1})=\Delta\Phi\left(\frac{\Delta}{\sigma_n'(\bm{x}_{n+1})}\right)+\sigma_n'(\bm{x}_{n+1})\phi\left( \frac{\Delta}{\sigma_n'(\bm{x}_{n+1})} \right)
\end{equation*}
Where $\Delta=T-\mu_{n}(\bm{x}_{n+1})$ and $T=\min\{\mu(\bm{x}_1),...,\mu(\bm{x}_n) \}$.
We now claim that $\bm{x}_{n+1}$ cannot be too close to any of the current evaluated points.
Otherwise, suppose  $\mu_n(\bm{x}_i)=T$ for $i\in U\subseteq \{1,2,3,...,n \}$.
If $\bm{x}_{n+1}$ tends to some point in $U$, then $\Delta\to 0$ and $\sigma_n'(\bm{x}_{n+1})\to 0$. By definition we know $ \text{EI}_T(\bm{x}_{n+1},\bm{P}_{n+1})\to 0$.
On the other hand, if $\bm{x}_{n+1}$ tends to some point $x_k$ in $\{1,2,3,...,n \}/ U$, then $\Delta\to T-\mu_n(\bm{x}_k)<0$ and $ \sigma_n'(\bm{x}_{n+1})\to 0$.
As a result, $ \Delta\Phi\left(\frac{\Delta}{\sigma_n'(\bm{x}_{n+1})}\right)\to ( T-\mu_n(\bm{x}_k) )\Phi(-\infty)=0$ and so $\text{EI}_T(\bm{x}_{n+1},\bm{P}_{n+1})\to 0. $ 
Thus, according to our stopping rule that the algorithm will stop when $\text{EI}_T <c$ for any candidate point for some value of $c$, we cannot select a $\bm{x}_{n+1}$ that is too close to any of the visited points.
From Theorem 1 of \cite{locatelli1997bayesian}, the algorithm will terminate within a finite number of steps. Lemma 1 of \cite{locatelli1997bayesian} implies that if we decrease the the value of $c$ to zero, the points we visited will be dense. 

As the evaluation points we visited are dense, there will exist a subsequence $\bm{x}_{n_k},\ k=1,2,...$ such that $\bm{x}_{n_k}\to_p \bm{x}_s$. Denote the estimated estimator at $n_k$-th iteration by $\widehat{\bm{x}}_{n_k}^*$. For $n_i < N < n_j$, as $n_i$ and $N$ go to infinity, we have
\begin{equation}
\begin{aligned}
\mu_N(\hat{\bm{x}}_N^*)\leq \mu_N(\hat{\bm{x}}_{n_i}^*) &\leq g(\hat{\bm{x}}_{n_i}^*) + o_p(1) = g(\bm{x}_s) + (g(\hat{\bm{x}}_{n_i}^*) -g(\bm{x}_s)) + o_p(1) = g(\bm{x}_s) + o_p(1)
\end{aligned}
\label{equ:inequality}
\end{equation}
The first inequality $\mu_N(\hat{\bm{x}}_N^*)\leq \mu_N(\hat{\bm{x}}_{n_i}^*)$ is due to the definition of $\hat{\bm{x}}_N^*$. The second inequality is due to the uniform convergence of  $|\mu_n(\bm{x})-g(\bm{x})|\to_p 0$ . The last equality is due to the continuous mapping theorem. With  (\ref{equ:inequality}), we have 
\begin{equation*}
g(\hat{\bm{x}}_N^*)  \leq \mu_N(\hat{\bm{x}}_N^*) + o_p(1) \leq g(\bm{x}_s) + o_p(1). 	
\end{equation*}

\textbf{Proof of Corollary \ref{lem: consistency of xN to xg}}
Now since $\bm{x}_s$ is unique, we can use the fact that $g$ is continuous at $\bm{x}_s$ to see 
\begin{equation*}
\inf_{|\bm{x}-\bm{x}_s|>\epsilon} g(\bm{x})>g(\bm{x}_s)
\end{equation*}
for all $\epsilon>0$.
And so
for any $\epsilon>0$, there exist a $\delta>0$ such that
\begin{equation*}
g(\bm{x})-g(\bm{x}_s)>\delta, \quad \mbox{for all} \quad |\bm{x}-\bm{x}_s|>\epsilon.
\end{equation*}

And so
\begin{equation*}
\mathbb{P}(|\hat{\bm{x}}_N^*-\bm{x}_s|>\epsilon)\leq \mathbb{P}(g(\hat{\bm{x}}_N^{*})-g(\bm{x}_s)>\delta)\to 0 , \quad \text{as }  N \to \infty.
\end{equation*}
The last step holds due to Theorem \ref{thm:main theorem}.

\textbf{Proof of Corollary \ref{corollary:2}}
In view of part $(ii)$ of Theorem \ref{theorem:1}, we can first choose $s_0$ such that, for any $s>s_0$,
\begin{equation}\label{eq: triangle1}
P(|\bm{x}_{s}-\bm{x}^*|\geq \frac{\epsilon}{2})<\frac{\epsilon}{2}.
\end{equation}
Now for any selected $s$, we use Corollary \ref{lem: consistency of xN to xg} to see that there exist $N_0$ such that, for any $N\geq N_0$,
\begin{equation}\label{eq: triangle2}
P(|\hat{\bm{x}}_{N}^*-\bm{x}_{s}|\geq \frac{\epsilon}{2})<\frac{\epsilon}{2}.
\end{equation}
Combining \eqref{eq: triangle1}-\eqref{eq: triangle2}, the following inequality is immediate:
\begin{equation*}
\mathbb{P}(|\hat{\bm{x}}_N^*-\bm{x}^*|\geq \epsilon)\leq P(|\hat{\bm{x}}_{N}^*-\bm{x}_{s}|\geq \frac{\epsilon}{2})+P(|\bm{x}_{s} - \bm{x}^*|\geq \frac{\epsilon}{2})<\epsilon. 
\end{equation*}

For the convergence of the optimal value, we first use part $(iii)$ of Theorem \ref{theorem:1} to see that there exist $s_0$ such that, for any $s>s_0$,
\begin{equation}\label{eq: triangle3}
P(|g(\bm{x}_{s})-f(\bm{x}^*, \bm{P}^c)|>\frac{\epsilon}{2})<\frac{\epsilon}{2}.
\end{equation}
Then for any fixed $s$, we use Theorem \ref{thm:main theorem} to see that there exists
$N_0$ such that, for any $N>N_0$,
\begin{equation}\label{eq: triangle4}
P(|g(\hat{\bm{x}}_N^*)-g(\bm{x}_{s})|>\frac{\epsilon}{2})<\frac{\epsilon}{2}.
\end{equation} 
Combining \eqref{eq: triangle3}-\eqref{eq: triangle4}, the following inequality is immediate:
\begin{equation*}
\resizebox{1\hsize}{!}{$
	\begin{aligned}
	P(|g(\hat{x}_{N}^*)-f(\bm{x}^*, \bm{P}^c)|\geq \epsilon)&<P(|g(\hat{\bm{x}}_N^*)-g(\bm{x}_{s})|>\frac{\epsilon}{2})  + P(|g(\bm{x}_{s})-f(\bm{x}^*, \bm{P}^c)|>\frac{\epsilon}{2}) <\epsilon.
	\end{aligned}
	$}
\end{equation*}

	\section{Test Functions} \label{appendix:function}
	\begin{table*}[!h]
		\centering
		\caption{Analytical expressions  of the test functions}
		\label{tab:info-test-func}
		{\small
			\begin{tabular}{|C{2cm}|C{12cm}|}
				\hline
				Function & Analytical expression \\
				\hline
				Griewank & 
				\begin{minipage}{12cm}
					\begin{equation*}
						f(x_1,x_2, u)=\left(\sum_{i=1}^2\frac{x_i^2}{4000} + \frac{u^2}{4000}-\cos(\frac{u}{\sqrt{3}})\prod_{i=1}^2\cos{\left(\frac{x_i}{\sqrt{i}}\right)} - 0.49\right)/0.48
					\end{equation*}
					$x_1,x_2\in[-50,50]$\\
					$u \in [-\infty, \infty]$\\  
					$u \sim 0.5 \text{lognorm} (10, 10) + 0.5 \text{lognorm}(20,5)$ 
					
				\end{minipage} \\
			   \hline

				StybTang & 
				\begin{minipage}{12cm}
					\begin{equation*}
						f(x_1,x_2, u)= \left(\frac{1}{2}(\sum_{i=1}^2(x_i^4- 16 x_i^2 + 5x_i) + u^4 - 16 u^2 + 5u )  - 398184\right)/17287676
					\end{equation*}
					$x_1,x_2\in[-5,5]$\\
					$u \in [-\infty, \infty]$\\
					$u \sim 0.5 \text{lognorm} (10, 20) + 0.5 \text{lognorm}(20,10)$ 
				\end{minipage} \\
				\hline
			\end{tabular}
		}
	\end{table*}

\end{document}